          \documentclass{article}
\usepackage{subfig}
\usepackage{amsmath}
\usepackage{amssymb}
\usepackage{graphicx}
\usepackage{dcolumn}
\usepackage{mathtools}
\usepackage{amsmath,amsfonts}
\usepackage{bm}
\usepackage{amsmath}
\usepackage{amssymb}
\usepackage{algorithm2e}
\usepackage{color}
\usepackage{algorithmic}
\usepackage{float}

\usepackage[paperwidth=7.5in, paperheight=10in, margin=.875in]{geometry}

            \newtheorem{thm}{Theorem}[section]

          \newtheorem{rem}[thm]{Remark}
          
\newcommand{\p}{\partial}
\newcommand{\f}{\frac}
\newcommand{\ds}{\displaystyle}
\newcommand{\E}{ {\mathbb{E}} }
\newcommand{\ola}{\overleftarrow}

\newcommand{\R}{\mathbb{R}}
\newcommand{\be}{\begin{equation}}
\newcommand{\ee}{\end{equation}}

 \date{}
           
\begin{document}

\title{Backward SDE Filter for Jump Diffusion Processes and Its Applications in Material Sciences  }

\author{ Richard Archibald  \thanks{  Computer Science and Mathematics Division, Oak Ridge National Laboratory,  Oak Ridge, Tennessee \ ({\tt archibaldrk@ornl.gov}).}
\and Feng Bao \thanks{  Department of Mathematics, The University of Tennessee at Chattanooga, Chattanooga, Tennessee, 37403 \ ({\tt feng-bao@utc.edu}).}
        \and  Peter Maksymovych \thanks{ Center for Nanophase Material Sciences, Oak Ridge National Laboratory, Oak Ridge, Tennessee  \ ({\tt maksymovychp@ornl.gov}).}
        }

\pagestyle{myheadings} 

\markboth{Optimal Filtering of Diffusion Processes}{Bao, Cao \& Han} \maketitle

\begin{abstract}
The connection between forward backward doubly stochastic differential equations and the optimal filtering problem is established without using the Zakai's equation.    The solutions of forward backward doubly stochastic differential equations are expressed in terms of conditional law of a partially observed Markov diffusion process.   It then follows that the adjoint time-inverse forward backward doubly stochastic differential equations governs the evolution of the unnormalized filtering density in the optimal filtering problem.
\end{abstract}


\markboth{Backward SDE Filter for Jump Diffusion Processes}{Archibald \& Bao \& Maksymovych} \maketitle

\begin{abstract}
We develop a novel numerical method for solving the nonlinear filtering problem of jump diffusion processes. The methodology is based on numerical approximation of backward stochastic differential equation systems driven by jump diffusion processes and we apply adaptive meshfree approximation to improve the efficiency of numerical algorithms. We then use the developed method to solve atom tracking problems in material science applications. Numerical experiments are carried out for both classic nonlinear filtering of jump diffusion processes and the application of nonlinear filtering problems in tracking atoms in material science problems.  
\end{abstract}

\textbf{keywords:} Nonlinear filtering problem, backward SDEs,  jump diffusion processes, material sciences


\section{Introduction}

The nonlinear filtering problem is one of the key missions in data assimilation, in which observations of a system are incorporated into the state of a numerical model of that system. Mathematically, the nonlinear filtering problem is to obtain, recursively in time, the best estimate of the state of unobservable stochastic dynamics $S = \{S_t: t \geq 0\}$, based on an associated observation process, $ M =\{M_t: t \geq 0\}$, whose values are a function of $S$ after corruption by noises. This suggests the optimal filtering problem of obtaining the conditional distribution of the state  $S_t$ from the observations up until time $t$, which achieves the best estimate of this distribution, in the squared error sense, based on the available observations. 

The nonlinear  filtering theory finds its applications in numerous scientific and engineering research areas, such as target tracking \cite{Djogatovi2014},  \cite{Kim2013}, signal processing  \cite{Hairer2011}, \cite{Little-Jones2013}, image processing \cite{Singh2013},  biology \cite{Dawson2015}, \cite{Lee2013}, \cite{Stepanov2015},  or mathematical finance \cite{Bensoussan2009},  \cite{Elliott2013}, \cite{Frey2012}.  
Some of the pioneer contributions to the development of nonlinear filters are due to Kushner \cite{MR0180407} and Stratonovich \cite{stratonovich1960}. Later, Zakai  \cite{zakai} introduced an
alternative approach to the computation of the nonlinear  filter by developing the so-called Zakai
equation, which is a stochastic partial differential equation (SPDE), and the best estimate of the nonlinear filter, i.e. the conditional distribution, is represented by the solution of the Zakai equation.  Although the Zakai's approach produces the ``\textit{exact}'' solution of the nonlinear filtering problem in theory, solving the SPDEs numerically can be extraordinarily difficult, especially when the state processes are in high dimensions \cite{BaoC2015}, \cite{Gobet-Zakai}, \cite{HU-Zakai}, \cite{zhang}. A more widely accepted method by practitioners to solve the nonlinear filtering problem is the sequential Monte Carlo approach, which is also known as the particle filter method \cite{Bao_KdV}, \cite{particle-filter-resample}, \cite{Chorin2013}, \cite{Doucet2008},  \cite{particle-filter}, \cite{Kang_2018}. The particle filter method uses a number of independent random variables, called particles, sampled directly from the state space to represent the prior probability, and updates the prior by including the new observation to get the posterior.  
This particle system is properly located, weighted and propagated recursively 
according to Bayes' theorem.  As a Monte Carlo approach, with sufficient large number
of samples the particle filter provides an accurate representation of the state probability density function (pdf)
as desired in the nonlinear filtering problem. Convergence of a particle filter to the optimal filter was shown under certain conditions \cite{MR1847785}, \cite{MR1895071}, \cite{HanX2015}. In addition to the Zakai's approach and Monte Carlo type approach, the authors have developed an alternative method, which solves the nonlinear filtering problem through a  forward backward doubly stochastic differential equations (BDSDEs) system. The theoretical basis of the BDSDEs approach is the fact that the BDSDEs system is equivalent to a parabolic type SPDE and the solution of that system is the conditional distribution of the state as required in the nonlinear filtering problem\cite{BaoC20142}, \cite{Bao_BDSDE_Filter}, \cite{Bao_Vasilis}. In this connection, it produces the exact solution of the nonlinear filtering problem, just like the Zakai's approach. In the meantime, as an stochastic \textit{ordinary} differential equation (SDE) approach, it also relies on stochastic sampling, just like the particle filter method. Therefore, the BDSDEs approach builds the bridge between the Zakai's approach and the Monte Carlo type approach.

In this paper, we consider a more general nonlinear filtering problem -- the nonlinear filtering problem for jump diffusion processes, in which the state process $S_t$ is a jump diffusion process and the state dynamic is perturbed by both traditional Gaussian noises and other kinds of L\`evy type noises. Different from classical nonlinear filtering problems, numerical methods to solve the nonlinear filtering problem for jump diffusion processes are not well developed. The existing methods for solving this type of problems focus on numerical approximation for its corresponding Zakai equation \cite{Jump-finance}, \cite{Popa-jump}, \cite{Duan}. However, due to the nonlocal behavior of the state dynamics as a  jump diffusion process, the corresponding Zakai equation contains fractional derivatives in spatial dimension, which is a stochastic partial-integral differential equation (PIDE).
The current numerical approaches to solve this type of stochastic PIDEs are extensions of existing methods for local partial differential equations, such like finite element, finite difference, etc..  It is well known that the partial-integral operator is a nonlocal operator, which may result in severe computational difficulties coming from the dramatic deterioration of the sparsity of the stiffness matrices required by the underlying linear systems.   Although the particle filter  framework could be extended to solve the nonlinear filtering problem for jump diffusion processes, as a typical drawback of Monte Carlo type approach, the particle filter method has very poor tail approximation for the state distribution especially when the state is perturbed by L\`evy type noises, which causes severe degeneracy problem for long term simulations \cite{Chen2003}, \cite{Doucet2008}.

In this work, we aim to develop an efficient numerical method to solve the nonlinear filtering problem for jump diffusion processes. The main theme of our approach is established on numerical approximations for a special kind of SDEs system, which is combined by a forward SDE and a backward SDE. We call this SDEs system a \textit{backward SDEs} system, and name the general backward SDEs approach for solving nonlinear filtering problems the \textit{Backward SDE filter}. Instead of using the equivalence relation between BDSDEs and the standard Zakai equation, in our novel approach for the nonlinear filtering problem for jump diffusion processes, we derive a backward SDEs system driven by jump diffusion processes to describe the time evolution of the state dynamics $S_t$. The theoretical validation of this effort is based on the probabilistic interpretation for PIDEs, which is discussed in \cite{BSDE-PIDE}, and our methodology in this paper simplifies the BDSDEs approach by avoiding simulations of the doubly stochastic term in BDSDEs. When receiving the observation measurements, we incorporate the data with the backward SDEs system by using the Bayes' theorem. In this connection, our method also produces a representation of the conditional state distribution as required in the nonlinear filtering problem and it maintains the accuracy property of the Zakai's approach. On the other hand, different from the Zakai's approach, the Backward SDE filter for jump-diffusion processes is still an \textit{SDE} based approach and it allows us to approximate the solution of nonlinear filtering problem on any selection of space points. Taking this advantage, we introduce a stochastic space points generation algorithm which generates meshfree space points adaptively according to the state dynamics and the observations. 
The central idea of this meshfree space points generation is to build dynamic random space points according to the approximate state distribution. 
Specifically, we generate a set of random samples from the initial distribution of the state and choose these samples as our random space points. Then, we propagate these space points through the dynamic system of the state. In this way, the space points move randomly according to the state model and are more concentrated around the state of the dynamic model in the nonlinear filtering problem. It is important to point out that the numerical approximation of the backward SDEs system for different selection of space points is independent. Therefore the Backward SDE filter has the same scalability as the Monte Carlo type approach. However, there is an essential difference with respect to the selection of random space points in the Backward SDE filter and the generation of random samples in the Monte Carlo type approaches. In the Backward SDE filter we approximate the value of state pdf on each space point instead of using the number of samples to describe the empirical distribution. In this way, the backward SDE filter requires much fewer space points compared to the sample-size in the Monte Carlo method.
With the numerical approximation of the solution of the backward SDE on the dynamic space points, we apply the Shepard's method which is an effective meshfree interpolation method to construct the approximation of the entire conditional pdf in the state space. In order to prevent space points degeneracy for long term simulations, we introduce a Markov Chain Monte Carlo based resampling method to update the space points according to the observation information.

Taking the accuracy and efficiency advantages of the Backward SDE filter, we then apply our developed method in solving problems in material science applications. Specifically, we study the application of Backward SDE filter in tracking atom trajectories moving on material surfaces. The so-called atom-tracking method provides a general approach to study diffusion of single atoms and molecules on a flat surface (typically a noble metal), and it also underlies the mechanism of single-atom manipulation  \cite{Heller_1994}, \cite{Hla_2014}, \cite{Morgenstern_2013}. In both cases an atomically sharp probe (made from a sharp metal needle) is scanned over the surface with sub-angstrom accuracy, typically producing the images of atomic-scale corrugation of the electronic density of states (in the case of scanning tunneling microscopy) or total electron density (in the case of atomic force microscopy \cite{Giessibl_2003}). The regimes of observation and manipulation are delineated in these methods of scanning probe microscopy by the strength of probe-atom (or probe-molecule) interaction, which in turn can be controlled by the degree of proximity between the probe and the observed atom  \cite{Peter_book}, \cite{Morgenstern_2013}. Many interesting examples of this method have been demonstrated -- from creation of artificial atomic structures (such as quantum corrals  \cite{Heller_1994}) to non-trivial diffusion mechanisms, including quantum tunneling of hydrogen. The tracking efficiency translates into both success and relative speed of manipulation, and is therefore of utmost importance. The frequency of both manipulation and diffusion of single atoms generally follow Poisson statistics, which makes the analysis of corresponding time-series a great fit to our Backward SDE filter. 

The major contribution of this paper is in twofolds: to develop an efficient Backward SDE filter for jump diffusion processes, and to introduce the potential application of the developed algorithms in material sciences. The rest of this paper is organized as following. In Section \ref{Prelim}, we introduce some preliminary theories for the nonlinear filtering problem and the backward SDEs system. In Section \ref{BSDE_Filter}, we develop the methodology and numerical algorithms for our Backward SDE filter. Numerical experiments are carried out in this section to examine the performance of our algorithms. In Section \ref{Material}, we discuss the application of the Backward SDE filter in material sciences and demonstrate the effectiveness of the Backward SDE filter in tracking atom trajectories on two different material surfaces.

\section{Preliminaries}\label{Prelim}
In this section, we introduce preliminary theories for the methodology of this paper. We first give a brief discussion to nonlinear filtering problems for jump diffusion processes. Then we introduce the backward SDEs system for jump diffusion processes and it's relation to integral-partial differential equations.
 
\subsection{Nonlinear filtering problems for jump diffusion processes}
Let $( \Omega, \mathcal{F},  \mathbb{P})$  be a probability space and denote $ \{\mathcal{F}_t\}_{0 \leq t \leq T}$ to be a filtration possesses right continuity, i.e. $\mathcal{F}_t = \mathcal{F}_{t+} $, and $\mathcal{F}_0$ is the $\sigma$-algrbra contains all the $\mathbb{P}$ zero measure zero sets. In this paper, the filtration $\{\mathcal{F}_t\}_{0 \leq t \leq T}$ is assumed to be generated by two mutually independent processes: a $d$-dimensional Brownian motion $W_t$ and a $d$-dimensional Poisson random measure $\mu(t, A)$ on $[0, T] \times \mathbb{E}$ where $\mathbb{E} = \mathbb{R}^d \backslash \{0\} $ is equipped with its Borel field $\mathcal{E}$, with compensator $\nu(dt, de) = dt \lambda(de)$, such that $\{\tilde{\mu} ([0,t] \times A) = (\mu - \nu) ([0,t] \times A) \}_{t \geq 0}$ is a martingale for all $A \in \mathcal{E}$ satisfying $\lambda(A) < \infty$ and $\lambda(de)$ is assumed to be a $\sigma$-finite measure on $(\E, \mathcal{E})$ satisfying
$$\int_{E} (1 \wedge |e|^2) \lambda (de) < + \infty,$$
where $| \cdot |$ denotes the standard Euclidean norm in Euclidean spaces.

In this paper, we consider the nonlinear filtering problem of jump diffusion processes in its state-space form on the probability space $( \Omega, \mathcal{F},  \mathbb{P})$ and introduce the following stochastic processes
\begin{subequations}\label{NLF}
\begin{align}
dS_t &= b(S_t) dt + \sigma_t dW_t + \int_E \beta_t(e) \tilde{\mu} (dt, de ),  \label{NLF:State} \\
M_t & = h(S_t, B_t), \label{NLF:Measurement}
\end{align}
\end{subequations}
where $b: \mathbb{R}^d \rightarrow \mathbb{R}^d$ and $h: \mathbb{R}^d \times \mathbb{R}^l \rightarrow \mathbb{R}^l$ are nonlinear functions, $\beta_t: \E \rightarrow \mathbb{R}^{d} $ is a $d$-dimensional process on $\E$, $W_t \in \mathbb{R}^d$ and $B_t \in \mathbb{R}^l$ are two independent Brownian motions, $\tilde{\mu}$ is the compensated Poisson random measure describes jumps in the model, $\sigma_t \in \mathbb{R}^{d \times d}$ is the coefficient matrix for Brownian motion $W_t$. The given initial value $S_0$ of \eqref{NLF:State} is independent of $W_t$, $\tilde{\mu}$ and $B_t$, and has a probability distribution with density function $Pr$. The stochastic process $\{S_t\}_{t\geq 0}$ defined by the stochastic differential equation (SDE) \eqref{NLF:State} describes the dynamics of a jump diffusion process, which is named the \textit{state process}, and $\{M_t\}_{t \geq 0}$ is the noise perturbed measurement of the state $\{S_t\}_{t \geq 0}$ with observation function $h(S_t, B_t)$.  When the coefficient $\beta \equiv 0$ in the state process \eqref{NLF:State}, the state-space form \eqref{NLF} gives a nonlinear filtering problem of standard Brownian motion driven diffusion processes.
In most applications, the measurement $M_t$ is received at discrete time and the the noise from measurement can be
assumed to be additive. In this way, the measurement is formulated in a discrete manner as
\begin{equation}\label{NLF:Discrete}
M_{t_n} = h(S_{t_n}) + R \dot{B}_{t_n}, \quad n = 1, 2, \cdots
\end{equation} 
where $R \in \mathbb{R}^l \times \mathbb{R}^l$ is the covariance matrix of the noise. 
The goal of the nonlinear filtering problem is to derive the least square estimate of a functional $\psi(S_t)$ of the state process $S_t$, given all the observation information of $M_t$ up to time $t$. In other words, we aim to find the optimal filter $\tilde{\psi}(S_t)$, such that
$$\tilde{\psi}(S_t) :=E[\psi(S_t) | \mathcal{M}_t] = \inf\{E[|\psi(S_t) - K_t|^2]; K_t \in \mathcal{K}_t\},$$ 
where $\mathcal{M}_t : = \sigma\{M_s, 0 \leq s \leq t\}$ is the $\sigma$-algebra containing all the observation information up to time $t$, and $\mathcal{K}_t$ is the space of all $\mathcal{M}_t$ measurable and square integrable random variables.  

\subsection{Backward SDEs driven by jump diffusion processes }

As preliminaries for our methodology, we introduce a system of backward SDEs on probability space $( \Omega, \mathcal{F},  \mathbb{P})$ driven by the state process $S_t$ defined in \eqref{NLF:State}.
Let $\mathcal{S}^2$ denote the set of $\mathcal{F}_t$-adapted c\`adl\`ag one dimensional processes such that for stochastic process $\{Y_t, 0 \leq t \leq T\} \in \mathcal{S}^2$, we have 
$$\|Y\|_{\mathcal{S}^2} = \| \sup_{0 \leq t \leq T}|Y| \|_{L^2(\Omega)} < \infty .$$
Let $L^2_W$ be the set of $\mathcal{F}_t$- progressively measurable $d$-dimensional processes such that for stochastic process $\{Z_t, 0 \leq t \leq T\} \in L^2_W$, we have 
$$\|Z\|_{\mathcal{S}^2} = \left( E \int_{0}^{T} | Z_t |^2 dt \right)^{1/2} < \infty .$$ 
In addition, we use $L^2_{\tilde{\mu}}$ to denote the set of mappings $U: \Omega \times [0, T] \times E \rightarrow \mathbb{R} $ which are $\mathcal{P} \otimes \mathcal{E}$ measurable, where $\mathcal{P}$ denotes the $\sigma$-algebra of $\mathcal{F}_t$-predictable subsets of $\Omega \times [0, T]$,  and such that 
$$\| U \|_{L^2(\tilde{\mu})} := \left( E \int_{0}^{T} \int_{E} U_t(e)^2 \lambda(de)  dt \right)^{1/2} < \infty .$$
Finally, we define $\mathcal{B}^2 = \mathcal{S}^2 \times L^2_W \times L^2_{\tilde{\mu}}$.

We consider the following forward backward stochastic differential equations system
\begin{subequations}\label{FBSDE:Feynman-Kac}
\begin{align}
\bar{X}_t =& \ \bar{X}_0 +  \int_{0}^{t} b(\bar{X}_s) ds + \int_{0}^{t}\sigma_s dW_s + \int_{0}^{t}\int_E \beta_s( e) \tilde{\mu} (ds, de ), \quad 0 \leq t \leq T,  \label{SDE:FK} \\
Y_t =& \psi(\bar{X}_T) + \int_{t}^{T} f_s(\bar{X}_s, Y_s, Z_s) ds - \int_{t}^{T} Z_s dW_s - \int_{t}^{T} \int_{E} U_s(e) \tilde{\mu} (ds, de), \quad 0 \leq t \leq T, \label{BSDE:FK}
\end{align}
\end{subequations}
where $W_t \in \mathbb{R}^d$ is a $d$-dimensional standard Brownian motion, $\tilde{\mu}$ is a compensated Poisson random measure, $b: \mathbb{R}^d \rightarrow \mathbb{R}^d$, $\sigma_t \in \mathbb{R}^{d \times d}$, $\beta_t \in \mathbb{R}^d$ and $f: [0, T] \times \mathbb{R}^d \times \mathbb{R} \times \mathbb{R}^d \rightarrow \mathbb{R}$ is a measurable function. For the convenience of presentation, we call the stochastic system \eqref{FBSDE:Feynman-Kac} a backward SDEs system. Equation \eqref{SDE:FK} in the above backward SDEs system is a forward SDE which describes the same diffusion of the state process \eqref{NLF:State} and equation \eqref{BSDE:FK} is a backward SDE with a given initial condition $Y_T = \psi(\bar{X}_T)$. The solution of the backward SDEs system \eqref{FBSDE:Feynman-Kac} is a quadruplet $(\bar{X}_t, Y_t, Z_t, U_t)$ with the triple $( Y_t, Z_t, U_t ) \in \mathcal{B}^2$. The existence and uniqueness of the solution quadruplet are provided by \cite{BSDE-PIDE}.   It is known that $Y_t$ and $Z_t$ are both functions of $\bar{X}_t$ (see \cite{Pardoux1990}) and 
we also known from \cite{BSDE-PIDE} that under the condition $\bar{X}_0 = x \in \mathbb{R}^d$,  
\begin{equation}
Z_t = \sigma_t \nabla Y_t, \hspace{1em} U_t = Y_t(\bar{X}_{t-} + \beta_t(e)) - Y_t(\bar{X}_{t-}),
\end{equation}
where $\nabla$ is the gradient and the minus sign in $X_{t-}$ is the left limit.

An important property of the backward SDEs system is its equivalence to the following parabolic integro-partial differential equation (PIDE)
\begin{equation}\label{PIDE}
\begin{aligned}
& - \f{\p u_t}{\p t} (x)= \mathcal{L}_t u_t(x) + f_t(x, u_t, \sigma_t \nabla u_t(x) ),  \quad (t, x) \in [0, T] \times \mathbb{R}^d,\\
& u_T(x) = \psi(x), \quad x \in \mathbb{R}^d,
\end{aligned}
\end{equation} 
where the second order integral-differential operator $\mathcal{L}$ is of the form $\mathcal{L}_t = A_t + K_t$, with
$$A_t\phi(x) = \f{1}{2} \sum_{i,j=1}^d\Big( (\sigma_t \sigma_t^{\ast} )_{ij}  \f{\p^2 \phi}{\p x^2} (x)\Big) + \sum_{i=1}^{d} (b_t)_i(x) \f{\p \phi}{\p x_i}(x), \quad \phi \in C^2(\mathbb{R}^d), $$
$$K_t\phi(x) = \int_{E} \big[ \phi\big(x + \beta_t(e)\big) - \phi(x) - \beta_t(e) \nabla \phi(x) \big]\lambda (de), \quad \phi \in C^2(\mathbb{R}^d).$$
The side condition of the above equation is given at time $T$ and the propagation direction is backward, i.e. from $T$ to $0$.
It has been shown in \cite{BSDE-PIDE} that the solution $Y$ of the backward SDEs system \eqref{FBSDE:Feynman-Kac} is the unique viscosity solution of \eqref{PIDE}, i.e. we have 
\begin{equation}\label{Y=u}
Y_t(x) = u_t(x),  \quad (t, x) \in [0, T] \times  \mathbb{R}^d.
\end{equation}


%

\section{Backward SDE filter for jump diffusion processes}\label{BSDE_Filter}
Now, we introduce the backward SDEs approach for the nonlinear filtering problem of jump diffusion processes. In subsection \ref{Methodology}, we first derive the methodology of solving the nonlinear filtering problem through a backward SDEs system and we name this methodology the \textit{Backward SDE filter for jump diffusion processes}. More generally, we call the backward SDEs type approach for the nonlinear filtering problem the \textit{Backward SDE filter}. In subsection \ref{Algorithms}, we develop an effective and efficient numerical scheme to solve the Backward SDE filter. Then, we present numerical examples to show the performance of our Backward SDE filter in solving the nonlinear filtering problem \eqref{NLF} in subsection \ref{Synthetic}.

\subsection{Methodology}\label{Methodology}

The framework of the Backward SDE filter is composed by two major tasks: (i) estimate the state process without observation information; (ii) incorporate observation with the state estimation.  The second task is typically achieved by applying Bayes' theorem, which we also adopt in the Backward SDE filter for jump diffusion processes. To achieve the first task, we introduce a backward SDEs system, which describes the probability density function (pdf) of the state, and solve the corresponding backward SDEs system. 


To proceed, we first consider the following special case of backward SDEs system \eqref{FBSDE:Feynman-Kac}
\begin{subequations}\label{FBSDE:Kolmogorov}
\begin{align}
S_t =& \ S_0 +  \int_{0}^{t} b(S_s) ds + \int_{0}^{t}\sigma_s dW_s + \int_{0}^{t}\int_E \beta_s( e) \tilde{\mu} (ds, de ), \quad 0 \leq t \leq T,  \label{SDE:Kolmogorov} \\
Y_t =& \psi(S_T) - \int_{t}^{T} Z_s dW_s - \int_{t}^{T} \int_{E} U_s(e) \tilde{\mu} (ds, de), \quad 0 \leq t \leq T. \label{BSDE:Kolmogorov}
\end{align}
\end{subequations}
where the forward SDE \eqref{SDE:Kolmogorov} is a standard jump diffusion process which has the same definition as the state process \eqref{NLF:State} in the nonlinear filtering problem.  Taking conditional expectation to both sides of \eqref{BSDE:Kolmogorov}, from the martingale property of $\int_{t}^{T} \cdot dW_t$ and $\int_{t}^{T} \int_{E} \cdot \tilde{\mu} (ds, de)$, we obtain
$$ E[Y_t] = E[\psi(S_T)], $$
which is a simplified Feynman-Kac formula.
From the equivalence relation between backward SDEs systems and PIDEs, we know that the backward SDEs system \eqref{FBSDE:Kolmogorov} is equivalent to the following PIDE
\begin{equation}\label{PIDE:Kolmogorov}
\begin{aligned}
 - \f{\p u_t}{\p t}(x) =& \f{1}{2} \sum_{i,j=1}^d\Big( (\sigma_t \sigma_t^{\ast} )_{ij}  \f{\p^2 u_t}{\p x^2} (x)\Big) + \sum_{i=1}^{d} (b_t)_i(x) \f{\p u_t}{\p x_i}(x) \\
& \qquad  + \int_{E} \big[ u_t\big(x + \beta_t(e)\big) - u_t(x) - \beta_t(e) \nabla u_t(x) \big]\lambda (de) \quad (t, x) \in [0, T] \times \mathbb{R}^d,\\
 u_T(x) = & \psi(x), \quad x \in \mathbb{R}^d,
\end{aligned}
\end{equation} 
which is also known as a Kolmogorov backward equation of jump diffusion processes. 
By applying integration by parts formula and Taylor expansion, one can prove that the Kolmogorov backward equation \eqref{PIDE:Kolmogorov} is adjoint to the following Fokker Planck type equation \cite{Duan, Duan_F-P}
\begin{equation}\label{PIDE:Fokker-Planck}
\begin{aligned}
\f{\p p_t}{\p t}(x) = & \sum_{i,j=1}^{d} \f{1}{2} (\sigma_t \sigma_t^{\ast})_{ij} \f{\p^2 p_t}{\p x^2}(x) - \sum_{i=1}^d \Big(  (b_t)_i(x)  \f{\p p_t}{\p x_i}(x) + \f{\p (b_t)_i}{\p x_i}(x) p_t(x)  \Big) \\
  & \qquad + \int_{E} \Big( p_t(x - \beta_t(e)) -  p_t(x)  + \beta_t(e) \nabla p_t\Big) \lambda(de), \quad (t, x) \in [0, T] \times \mathbb{R}^d, \\
 p_0(x) =& Pr(x), 
\end{aligned}
\end{equation}
where the initial condition $p_0$ is the given pdf $Pr$ of the initial state random variable $S_0$, and the propagation direction of the  equation is forward, i.e. from $0$ to $T$. Actually, it is well known that the solution $p_t$ of the Fokker-Planck type equation \eqref{PIDE:Fokker-Planck} describes the probability evolution of the stochastic process $S_t$ defined by \eqref{NLF:State}, i.e. $p_t$ is the pdf for $S_t$. 

From the equivalence condition between backward SDEs and PIDEs, we know that there's a backward SDEs system corresponding to the Fokker-Planck type equation \eqref{PIDE:Fokker-Planck}. Since the Fokker-Planck type equation \eqref{PIDE:Fokker-Planck} has opposite time propagation direction to the Kolmogorov backward equation \eqref{PIDE:Kolmogorov}, the backward SDEs system equivalent to \eqref{PIDE:Fokker-Planck} also has opposite propagation direction to the backward SDEs system \eqref{FBSDE:Kolmogorov}. From the general equivalence relation between \eqref{FBSDE:Feynman-Kac} and \eqref{PIDE}, and the expression of the Fokker-Planck type equation \eqref{PIDE:Fokker-Planck}, we obtain that the following backward SDEs system is equivalent to equation \eqref{PIDE:Fokker-Planck}
\begin{subequations}\label{FBSDE:Fokker-Planck}
\begin{align}
X_0 =& \ X_t -  \int_{0}^{t} b( X_s) ds - \int_{0}^{t}\sigma_s d\ola{W}_s - \int_{0}^{t}\int_E \beta_s( e) \tilde{\mu} (ds, de ), \quad 0 \leq t \leq T,  \label{SDE:Fokker-Planck} \\
P_t =& Pr(X_0) - \int_{0}^{t} \sum_{i=1}^d \f{\p (b_s)_i}{\p x_i}( X_s) P_s ds - \int_{0}^{t} Q_s d\ola{W}_s - \int_{0}^{t} \int_{E} V_s(e) \tilde{\mu} (ds, de), \quad 0 \leq t \leq T, \label{BSDE:Fokker-Planck}
\end{align}
\end{subequations}
where $\int_0^t \cdot d \ola{W}_t$ is an It\^o integral integrated backwardly, i.e. from $t$ to $0$, adapted to a backward filtration of the Brownian motion $W_t$ and is named \textit{backward It\^o integral} (see \cite{Two_sided} for details).
With any given state $X_t$ as the side condition and the backward It\^o integral in the equation, \eqref{SDE:Fokker-Planck} can be considered as an SDE with backward propagation direction.  Similarly, we can observe that the second equation in \eqref{FBSDE:Fokker-Planck} is also a time inverse backward SDE with initial condition $P_0 = Pr$ and has the same structure to \eqref{BSDE:FK}.
As a result, \eqref{FBSDE:Fokker-Planck} is a time inverse backward SDEs system with solution quadruplet $(X_t, P_t, Q_t, V_t)$. 
From the equivalence relation \eqref{Y=u}, we know that the solution $P_t$ of the backward SDEs system \eqref{FBSDE:Fokker-Planck} is equivalent to the solution $p_t$ of the Fokker-Planck type equation \eqref{PIDE:Fokker-Planck}, i.e. 
\begin{equation}\label{P_t=p_t}
P_t(x) = p_t(x)
\end{equation}
Therefore, the numerical approximation of the solution $P_t$ in \eqref{FBSDE:Fokker-Planck} is also the  approximation for the pdf of state $S_t$. 

In this connection, we introduce the methodology of the Backward SDE filter, which solves for the conditional pdf of the state process in two steps: a \textit{Prediction Step} and an \textit{Update Step}. In the Prediction Step, we solve the backward SDE system \eqref{FBSDE:Fokker-Planck} numerically, and use the approximate solution as the predicted pdf for the state without using observation information. In the Update Step, we incorporate the observed measurement with the predicted pdf by using the Bayes' theorem.
In what follows, we provide our efficient numerical algorithms for the Backward SDE filter.

\subsection{Numerical algorithms}\label{Algorithms}

To introduce a recursive discretized numerical scheme, we consider a temporal partition over the time period $[0, T]$ as $\mathcal{R}_t := \{t_n | 0 = t_0 < t_1 < t_2 < \cdots < t_{N_T - 1} < t_{N_T} = T\}$ and let $\Delta t_n = t_{n+1} - t_n$, $\Delta W_{t_n} = W_{t_{n+1}} - W_{t_n}$. We also assume that the observation measurements are received at time $t_n$, $n = 1, 2, \dots, N_T$. 

\subsubsection*{\underline{\textbf{Backward SDE filter Framework}}}

\textit{Prediction Step.} \ 
The major task in the Prediction Step is to solve the backward SDEs system \eqref{FBSDE:Fokker-Planck} on time interval $[t_n, t_{n+1}]$, $n = 0, 1, \dots, N_T - 1$, i.e.
\begin{subequations}\label{FBSDE:t_n}
\begin{align}
X_{t_n} =& \ X_{t_{n+1}} -  \int_{t_n}^{t_{n+1}} b( X_s) ds - \int_{t_n}^{t_{n+1}}\sigma_s d\ola{W}_s - \int_{0}^{t}\int_E \beta_s( e) \tilde{\mu} (ds, de ),   \label{SDE:t_n} \\
P_{t_{n+1}} =& P_{t_n} - \int_{t_n}^{t_{n+1}} \sum_{i=1}^d \f{\p (b_s)_i}{\p x_i}( X_s) P_s ds - \int_{t_n}^{t_{n+1}} Q_s d\ola{W}_s - \int_{t_n}^{t_{n+1}} \int_{E} V_s(e) \tilde{\mu} (ds, de), \label{BSDE:t_n}
\end{align}
\end{subequations}
assuming $P_{t_n}$ is known. It's worthy to point out that in the nonlinear filtering problem, the initial condition $P_{t_n}$ for  the above backward SDEs system is chosen to be the conditional pdf of the state $S_{t_n}$ given observation information $\mathcal{M}_{t_{n}}$, i.e. $P_{t_{n}} = p(S_{t_n} | \mathcal{M}_{t_{n}})$.
In what follows, we first derive a temporal discretization scheme for \eqref{FBSDE:t_n}, and then discuss the spatial discretization later in this subsection.

Consider the equation \eqref{SDE:t_n}. We use Euler-Maruyama scheme to discretize the integrals and approximate the solution $X_{t_n}$ by
\begin{equation}\label{Semi_X}
X_{t_n} =  X_{t_{n+1}} - b( X_{t_{n+1}}) \Delta t_n - \sigma_{t_{n+1}} \Delta W_{t_n} - \int_E \beta_{t_{n+1}}( e) \tilde{\mu} (\Delta t_n, de ) + R_X^n,
\end{equation}
where $R_X^n$ is the approximation error of the equation.

To derive a numerical scheme for the backward SDE in \eqref{FBSDE:t_n}, we take conditional expectation $E_{n+1}[\cdot]$ on both sides of  \eqref{BSDE:t_n}, where $E_{n+1}[\cdot] : = E[\cdot \big| X_{t_{n+1}}]$. It follows from the facts 
$$E_{n+1}\left[\int_{t_n}^{t_{n+1}} Q_s d\ola{W}_s\right] = 0,$$
$$E_{n+1}\left[\int_{t_n}^{t_{n+1}} \int_{E} V_s(e) \tilde{\mu} (ds, de) \right] = 0$$
and $E_{n+1}[P_{t_{n+1}}] = P_{t_{n+1}}$
that the equation \eqref{BSDE:t_n} becomes
 \begin{equation}\label{Exp_P}
P_{t_{n+1}} = E_{n+1}[P_{t_n}] - \int_{t_n}^{t_{n+1}}E_{n+1}\Bigg[\sum_{i=1}^d \f{\p (b_s)_i}{\p x_i}( X_s) P_s \Bigg]ds.
\end{equation}
Then, we use the left point formula to discretize the deterministic integral in \eqref{Exp_P} to get
 \begin{equation}\label{Semi_P}
P_{t_{n+1}} = E_{n+1}[P_{t_n}] - E_{n+1}\Bigg[\sum_{i=1}^d \f{\p (b_{t_n})_i}{\p x_i}( X_{t_n}) P_{t_n} \Bigg] \Delta t_n + R_P^n,
\end{equation}
where $R_P^n$ is the approximation error of the equation.

Next, we drop the errors terms $R_X^n$ and $R_P^n$ in approximations \eqref{Semi_X} and \eqref{Semi_P}, respectively, and obtain the numerical schemes for solving $X$ and $P$ as following
\begin{subequations}\label{Scheme}
\begin{align}
&X_n = X_{n+1} - b(X_{n+1}) \Delta t_n - \sigma_{t_{n+1}} \Delta W_{t_n} - \int_E \beta_{t_{n+1}}( e) \tilde{\mu} (\Delta t_n, de ), \label{Scheme:X} \\
&\tilde{P}_{n+1} =  E_{n+1}[P_{n}] - E_{n+1}\Bigg[\sum_{i=1}^d \f{\p (b_{t_n})_i}{\p x_i}( X_{n}) P_{n} \Bigg] \Delta t_n, \label{Scheme:P}  
\end{align}
\end{subequations}
where $X_n$ is an approximation for solution $X_{t_n}$ and $\tilde{P}_{n+1}$ is an approximation for solution $P_{t_{n+1}}$.   From the equivalence relation \eqref{P_t=p_t}, we know that the numerical solution $\tilde{P}_{n+1}$ is an approximation for conditional pdf of the state $S_{t}$ at time instant $t_{n+1}$ before receiving measurement data, i.e. $\tilde{P}_{n+1} \approx p(S_{n+1} \big| \mathcal{M}_{t_n})$.

\begin{rem}
Although the solution of the backward SDE system \eqref{FBSDE:t_n} is a quadruplet $(X, P, Q, V)$, we do not need numerical approximations for solutions $Q$ and $V$ in neither the nonlinear filtering applications, nor the numerical scheme \eqref{Scheme}. Therefore, in this approach, we do not discuss numerical schemes for $Q$ and $V$. 
\end{rem}

\textit{Update Step.} \ 
To derive an approximation for the conditional pdf $p(S_{t_{n+1}}\big| \mathcal{M}_{t_{n+1}})$ and to incorporate the measurement data at time $t_{n+1}$, we apply Bayes' theorem to combine the estimate pdf $\tilde{P}_{n+1}$ obtained in the Prediction Step with the data $M_{t_{n+1}}$. Specifically, we let
\begin{equation}\label{Bayes}
\Pi_{n+1}(x) := \ds \f{\tilde{P}_{n+1}(x)\Psi_{n+1}(x)}{C_{n+1}},
\end{equation}
where 
$\Psi_{n+1}(x) : = \exp\left( -\f{1}{2 R} \| M_{t_{n+1}} - \psi(x) \|^2 \right)$ is proportional to the Gaussian likelihood function, $C_{n+1}$ is a normalization factor. Apparently,  $\Pi_{n+1}$ is an approximation for the conditional pdf of the state given observation information $M_{t_{n+1}}$, i.e. $\Pi_{n+1} \approx p(S_{t_{n+1}} \big| \mathcal{M}_{t_{n+1}})$, and we let $P_{n+1} = \Pi_{n+1}$.

\vspace{1em}

With the Prediction Step and the Update Step introduced above, we establish the basic framework of our Backward SDE filter: At each recursive time stage $t_{n} \rightarrow t_{n+1}$, $n = 0, 1, 2, \dots, N_{T}-1$, we let $P_{n} = \Pi_{n}$ be the initial condition of the time inverse backward SDE system and use numerical scheme \eqref{Scheme} to calculate predicted state pdf $\tilde{P}_{n+1}$. Then, we update the state pdf through \eqref{Bayes} to get the approximate conditional pdf $\Pi_{n+1}$ for the state $S_{t_{n+1}}$.

\subsubsection*{\underline{\textbf{Adaptive Meshfree Approximations.}}}

The numerical scheme \eqref{Scheme} for the backward SDE system \eqref{FBSDE:t_n} can be considered as a temporal discretization scheme. In order to provide an effective numerical approximation for the solution $P_{t_{n+1}}$ as a function of the random variable $X_{t_{n+1}}$, we need to approximate the conditional expectation $E_{n+1}[\cdot]$, which is a functional of $X_{t_{n+1}}$. This could also be considered as a spatial discretization method. Since $X_{t_{n+1}}$ is a continuous random variable in the state space $\R^{d}$, it's impossible to approximate $\E_{n+1}[\cdot]$ on all possible values of $X_{t_{n+1}}$ and a representation of $X_{t_{n+1}}$ is required. In this work, we choose a set of space points, denoted by $\mathcal{D}_{n+1}:= \{x_{n+1}^1, \dots, x_{n+1}^N\} \in \R^d$, to be a representation of $X_{t_{n+1}}$, and the corresponding conditional expectation values, i.e. $\{ E[\cdot | X_{t_{n+1}} = x^i_{n+1}] \}_{x^i_{n+1} \in \mathcal{D}_{n+1}}$, is the approximation of the conditional expectation $E_{n+1}[\cdot]$.  Although the standard tensor product mesh is a straightforward option for the representation of the random variable $X_{t_{n+1}}$, the numerical approximation for $E_{n+1}[\cdot]$ over a tensor product mesh is not feasible for two reasons. First of all, the tensor product mesh suffers the so-called ``curse of dimensionality'' -- the number of grid points increases exponentially as the dimension $d$ increases. Secondly, the pdf of the random variable $X_{t_{n+1}}$ has unbounded support and the domain of the tensor product mesh needs to be sufficiently large to cover the true target state. 

To address the aforementioned difficulties of tensor product mesh, we use a stochastic meshfree construction of $\mathcal{D}_{n+1}$, which could also be considered as an adaptive space points generation method. To proceed, we first generate a set of $N$ random samples, denoted by $\{\xi^i\}_{i=1}^{N}$, from the pdf $Pr$ of the initial state $S_0$. Apparently, the random samples $\{\xi^i\}_{i=1}^{N}$ are more concentrated in the high density region of $Pr$ and we let $\mathcal{D}_0 := \{x_0^i\} = \{\xi^i\}$, i.e. $x_0^i = \xi^i$ for $i=1, 2, \dots, N$. 
For time step $n =0, 1, 2, \dots, N_T-1$, we propagate space points $ \{x_{n}^i\}$ to $\{x_{n+1}^i\}$ through the state dynamic \eqref{NLF:State}, i.e.
\begin{equation}\label{State:Points}
x_{n+1}^i = x_n^i + b(x_n^i) \Delta t_n + \sigma_{t_{n}} w_{t_n}^i + L_{t_n}^i(\beta, \tilde{\mu}), \qquad i = 1, 2, \dots, N, 
\end{equation}
to get our space points set $\mathcal{D}_{n+1}: = \{x_{n+1}^i\}$ at time stage $t_{n+1}$, where $w_{t_n}^i$ is the $i$-th random sample according to the normal distribution $N(0, \Delta t_n)$, and we denote $L_{t_n}^i(\beta, \tilde{\mu})$ as the numerical approximation for the L\`evy term $\int_E \beta_{t_{n}}( e) \tilde{\mu} (\Delta t_n, de )$ corresponding to the $i$-th sample. There are different numerical approximation schemes for $L_{t_n}^i(\beta, \tilde{\mu})$ based on the choice of L\`evy characteristic function in the problem. The discussion of numerical approximations for L\`evy processes is out of scope of this paper and we refer to  \cite{Kloeden_Jump, Platen_Jump, Karniadakis_Temp} for details.  
We can see from our construction of space points that the points in $\mathcal{D}_{n+1}$, $n=0, 1, 2, \dots, N_T-1$, move dynamically according to the state model \eqref{NLF:State} in a stochastic manner and the points are more concentrated in the high probability density region of the state $S_{t_{n+1}}$.

With a set of points $\mathcal{D}_{n+1}$, we approximate the conditional expectation $E_{n+1}[\cdot]$ on each point in $\mathcal{D}_{n+1}$, i.e. compute $\{E_{n+1}^{x^i}[\cdot]\}_{i=1}^N$, where $E_{n+1}^{x^i}[\cdot]:= E[\cdot \big| X_{t_{n+1}}] \big|_{X_{t_{n+1}} = x_{n+1}^i}$ is the value of conditional expectation on the space point $x_{n+1}^i$.
In this research, we use the Monte Carlo method to approximate conditional expectations \cite{MCMC-PF, Monte-Carlo}.  To be specific, for each given space point $x_{n+1}^i$, we approximate the expectation $E_{n+1}^{x^i}[P_n]$ on the right hand side of \eqref{Scheme:P} as
\begin{equation}\label{Monte-Carlo}
E_{n+1}^{x^i}[P_n] \approx \f{1}{M} \sum_{m=1}^{M} \hat{P}_n(x_n^{i, m}), 
\end{equation} 
where $M$ is the number of samples in the Monte Carlo simulation, and the space point $x_n^{i,m}$ is obtained by solving the forward SDE in the backward SDE system, which is given by equation \eqref{Scheme:X} as
$x_n^{i, m} = x_{n+1}^i - b(x_{n+1}^i) \Delta t_n + \sigma_{t_{n+1}} \bar{w}_{t_n}^{m} - \bar{L}_{t_{n+1}}^m(\beta, \tilde{\mu})$, where $\bar{w}_{t_n}^{m} $ is the $m$-th sample of the distribution $N(0, \Delta t_n)$ independent from $\{w_{t_n}^i\}$ and $\bar{L}_{t_{n+1}}^m$ is the $m$-th sample of the numerical approximation for $\int_E \beta_{t_{n+1}}( e) \tilde{\mu} (\Delta t_n, de )$, which is also independent from $\{L_{t_n}^i\}$. The approximation term $\hat{P}_n(x_n^{i, m})$ in equation \eqref{Monte-Carlo} is an interpolatory approximation of $P_n$ at $x_n^{i, m}$ based on values $\{P_n(x_n^{i})\}_{x_n^i \in \mathcal{D}_n}$ with the scheme
$$\hat{P}_n(x_n^{i, m}) = \sum_{x_n^i \in \mathcal{D}_n} P_n(x_n^{i}) \Gamma^i (x_n^{i, m}),  $$
where $\{\Gamma^i\}_{i=1}^{N}$ is a set of basis functions. The Monte Carlo approximation for the expectation $E_{n+1}\Big[\sum_{i=1}^d \f{\p (b_{t_n})_i}{\p x_i}( X_{n}) P_{n} \Big]$ in \eqref{Scheme:P} is followed directly by a scheme similar to \eqref{Monte-Carlo}.

In order to approximate $\hat{P}_n$ given values of $P_n$ on randomly generated meshfree space points $\mathcal{D}_n$, an effective interpolation method is needed. However, the standard polynomial interpolation methods are not applicable to approximate $\hat{P}_n$ due to uncontrollable approximation errors\cite{Walter-Error}. To overcome this challenge, we use \textit{radial basis approximation} to construct the interpolant $\hat{P}_n$. Specifically, in this work we apply Shepard's method\cite{RBF} as our radial basis approximation to compute $\hat{P}_n(x)$ for any point $x \in \mathbb{R}^d$. The Shepard's method is also known as the ``Inverse Distance Weighting'' method. It uses the weighted average of values $\{P_n(x_n^{i})\}_{x_n^i \in \mathcal{D}_n}$ based on the distance between $x$ and $\{x_n^i\}$ to construct the interpolant. For a given space point $x \in \mathbb{R}^d$, we reorder the points in $\mathcal{D}_n$ by the distance to $x$ from short to long to get a sequence $\{x_n(j)\}_{j=1}^N$ such that $\| x_n(j) - x \| \leq \| x_n(k) - x \|$ if $j < k$, where $\| \cdot \|$ is the Euclidean norm in $\mathbb{R}^d$. Then, we choose a proper integer $J \leq  N$ and use the weighted average of the first $J$ values in $\{  P_n(x_n(j))\}_{j=1}^N$ to approximate $\hat{P}_n(x)$, i.e. 
\begin{equation*}\label{Shepard}
\hat{P}_n(x) = \sum_{j=1}^J P_n(x_n^{i}) \bar{\Gamma}^j (x),
\end{equation*}
where $\bar{\Gamma}^j (x)$ is the inverse distance weight and is defined by
$$\bar{\Gamma}^j (x) := \f{\| x_n(j) - x \|}{\sum_{j=1}^{J} \| x_n(j) - x \|},  \quad x_n(j) \in \{x_n(j)\}_{j=1}^N.$$
Note that $\sum_{j=1}^{J} \bar{\Gamma}^j (x) = 1$. 
As a result of our stochastic space points generation method and the meshfree approximation in Monte Carlo simulations,  we obtain numerical approximations for conditional expectations, and therefore numerical approximations for the state pdf on adaptively selected meshfree space points.  
 \vspace{0.75em}

It is worthy to point out that there are essential differences between the adaptive meshfree approximation in the Backward SDE filter and the empirical approximation in the Monte Carlo type approaches. Although both the stochastic space points generation in the Backward SDE filter and the traditional Monte Carlo method generate random samples from a given probability distribution, the Backward SDE filter could provide more accurate approximation for the solution of nonlinear filtering problems with fewer sample points for the following reasons. First of all, the Backward SDE filter approximates the value of the state pdf, which is the solution $P_t$ in the backward SDE system, at each space point instead of using the number of samples to describe an empirical distribution in the Monte Carlo type approach. In this way, the Backward SDE filter requires much fewer space points compare to the sample-size of the Monte Carlo method. Secondly, in the Backward SDE filter we use meshfree interpolation to construct a smooth approximation of the state pdf. This is unlike the empirical distribution obtained in the Monte Carlo method, which is equivalent to a piece-wise constant approximation. Therefore the approximation that we obtain in the Backward SDE filter is smoother and more accurate than the empirical distribution simulation in the Monte Carlo method. 
In addition, for nonlinear filtering problems of jump diffusion processes, the state distribution is more likely to have heavy tails due to the nonlocal behavior of the state process. In this case, people need sufficient large number of samples in the Monte Carlo method to describe the heavy tails in a relatively large region. On the other hand, the adaptive meshfree approximation in the Backward SDE filter would provide smooth and more accurate tail distributions with much fewer space points in the tails. 

To demonstrate the performance of the adaptive meshfree approximation in describing heavy-tailed distributions, we approximate an $\alpha$ stable distribution by using adaptive meshfree approximation and classic Monte Carlo method. 
To proceed, we define an $\alpha$ stable distribution $\Phi (\alpha, \gamma, \beta, \delta)$ of a random variable $X$ with characteristic function given by
$$
\begin{aligned}
E[\exp(i t X)] = \exp\Big( -\gamma^{\alpha} | t |^{\alpha} \Big[1 + i \beta  \text{sign}(t) \tan \f{\pi \alpha}{2} \big( (\gamma | t |)^{1-\alpha} - 1 \big)  \Big] + i \delta t \Big), \quad \alpha \neq 1,
\end{aligned}
$$ 
where the  parameters are chosen to be $\alpha = 0.75$, $\gamma = 1$, $\beta = 0$ and $\delta = 0$.
\begin{figure}[h!]
\begin{center}
\subfloat[ ]{\includegraphics[scale = 0.45]{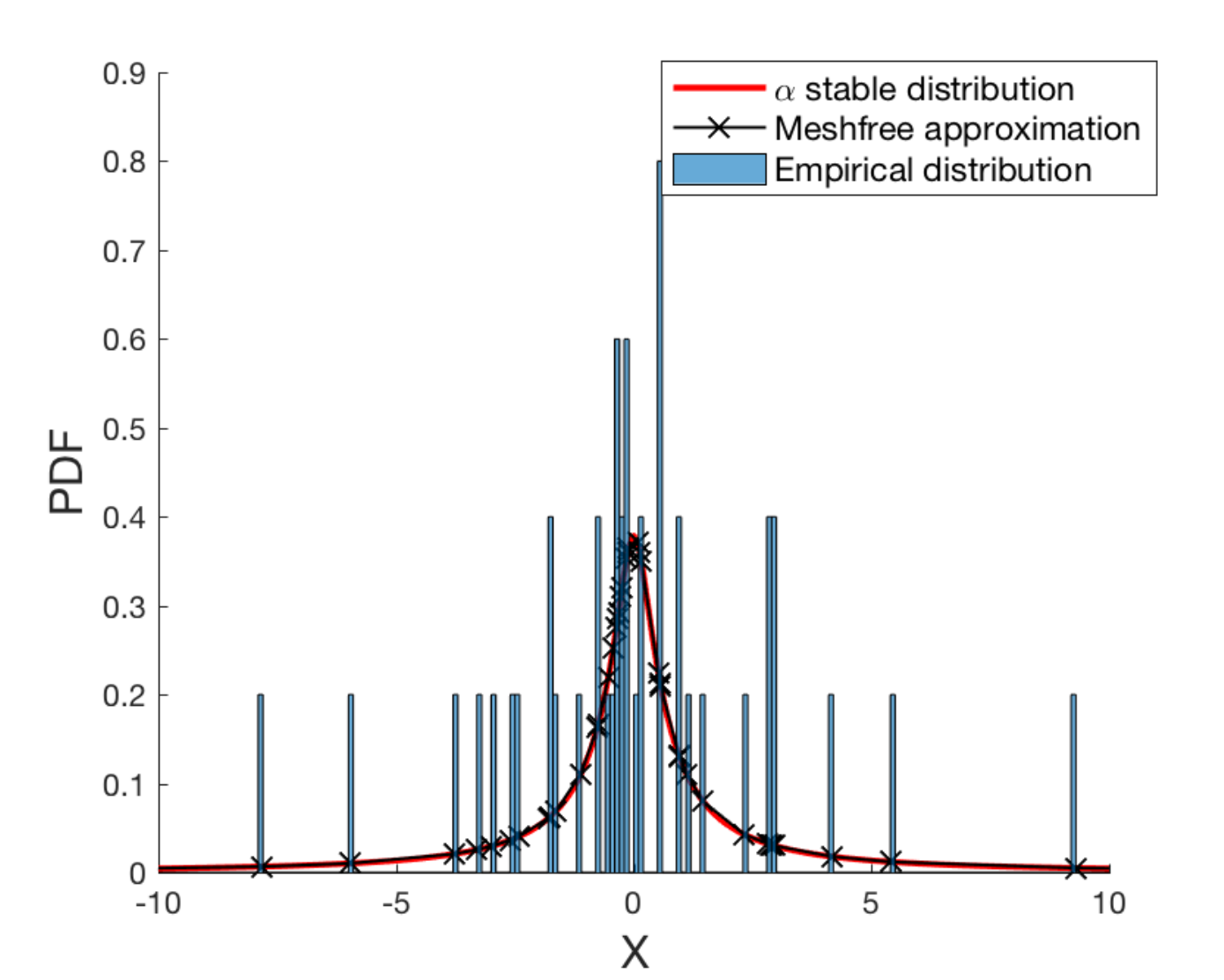} } \hspace{0.5em}
\subfloat[ ]{\includegraphics[scale = 0.45]{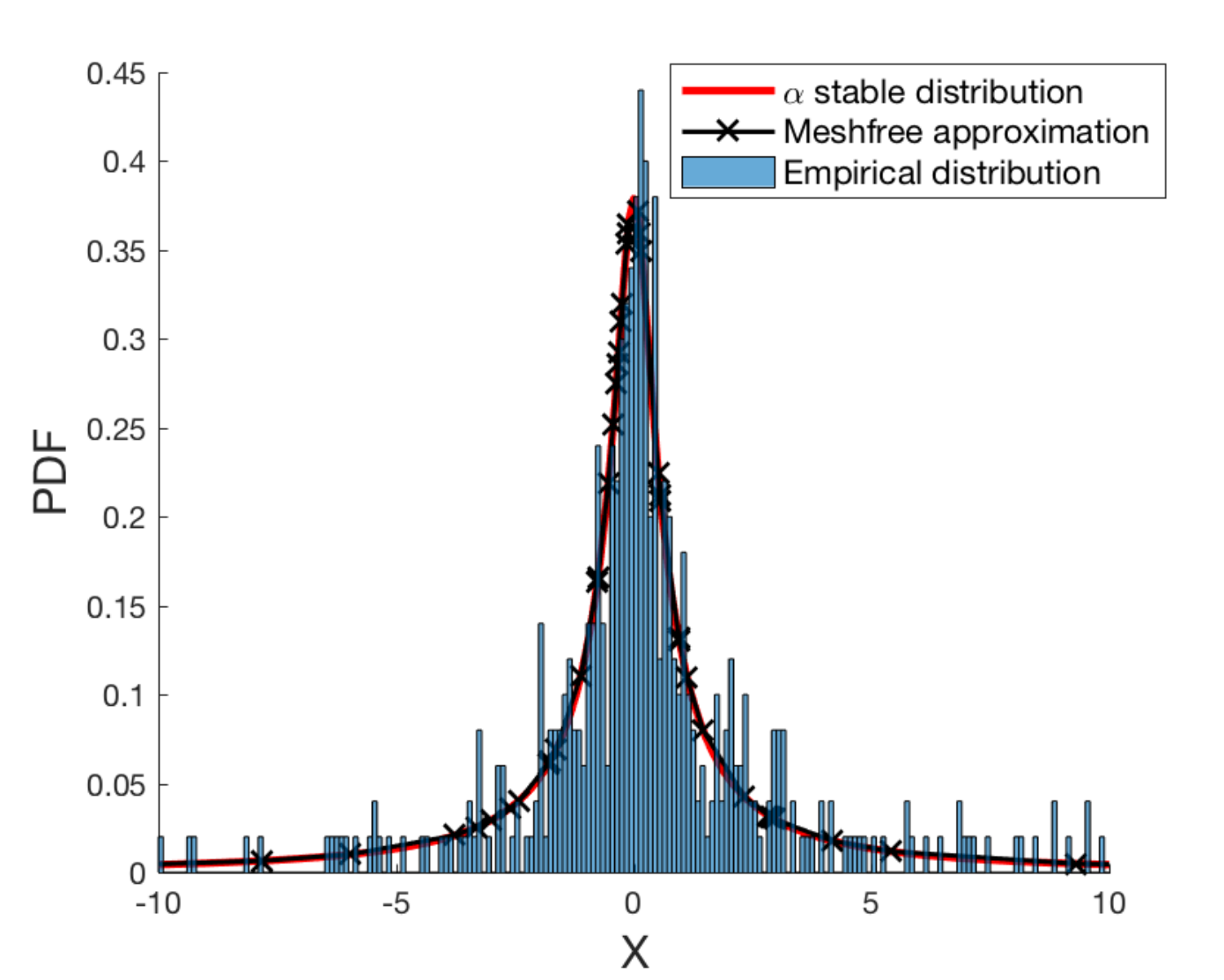} }\hspace{0.5em}
\subfloat[ ]{\includegraphics[scale = 0.45]{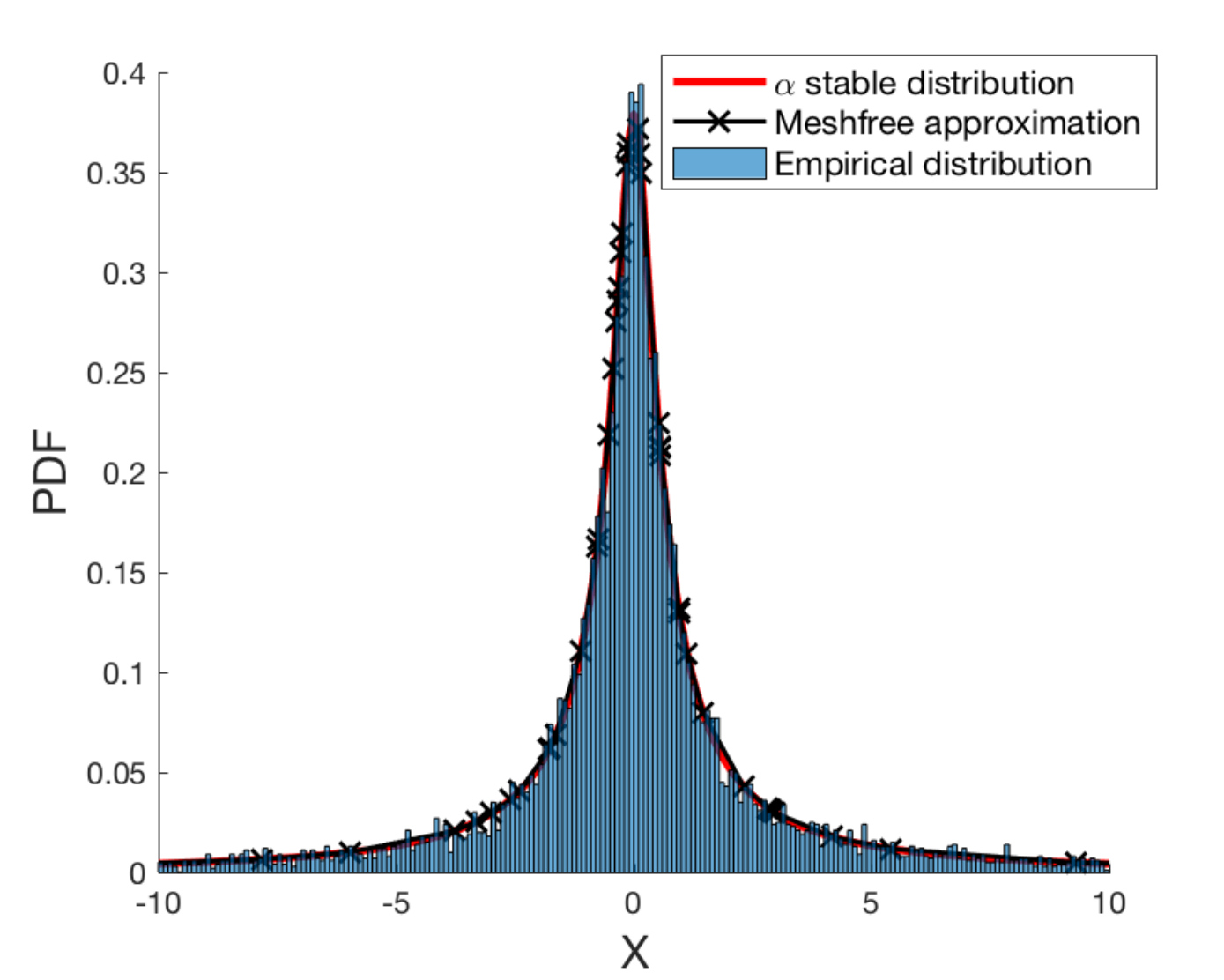} }
\end{center}
\vspace{-1em}
\caption{Comparison of adaptive meshfree approximation with the empirical distribution of Monte Carlo method in approximating the $\alpha$ stable distribution $\Phi$: (a) empirical distribution with $50$ samples; (b) empirical distribution with $500$ samples; (c) empirical distribution with $10,000$ samples. The adaptive meshfree approximation uses $50$ space points. }\label{Demo_distribution}
\end{figure}
In Fig. \ref{Demo_distribution}, we compare the adaptive meshfree approximation with the empirical distribution of Monte Carlo method in approximating the $\alpha$ stable distribution $\Phi$. The original $\alpha$ stable distribution is represented by the red curve in each subplot. For the adaptive meshfree approximation, we approximate the distribution with $50$ space points and use the black curve and black cross marks to describe the approximated distribution and the adaptive space points, respectively. The histogram in each subplot is the empirical distribution of Monte Carlo method and the number of samples increases from $50$ to $500$, and then to $10,000$ in subplots (a), (b) and (c), respectively. Since the adaptive space points in the meshfree approximation are chosen randomly according to the target distribution, in this demonstration we use the same $50$ samples as the Monte Carlo simulation presented in subplot (a). We can see from the figure that as the sample-size increases, the Monte Carlo simulation becomes smoother and more accurate, and the empirical distribution with $10,000$ samples is comparable to adaptive meshfree approximation with $50$ space points.  Moreover, we can see from subplots (a) and (b) that with smaller sample-size in Monte Carlo simulations, very few samples lie in the tails of distribution which makes the tail description very unreliable. On the other hand, although only $50$ space points are used  to approximate the distribution in the adaptive meshfree approximation and limited space points lie in the tail regions, the interpolatory approximation still makes the tail distribution smooth and accurate.

\subsubsection*{\underline{\textbf{Markov Chain Monte Carlo Resampling.}}}

In the stochastic space points generation method, the space points $\mathcal{D}_n$ move according to the state model \eqref{NLF:State}, which is a diffusion process. Therefore, the space points cloud diverge for long term simulations and less space points are located in the high probability density region of the state pdf after several simulation steps. This would make the space points very sparse in the state space and the probability density tends to concentrate on a few points, which dramatically reduces the meshfree interpolation accuracy.

To avoid the divergence problem of stochastic space points, inspired by the resampling procedure in the particle filter method \cite{MCMC-PF}, \cite{particle-filter}, we introduce a Markov Chain Monte Carlo (MCMC) method based resampling step. MCMC method is a class of algorithms for sampling from a probability distribution based on constructing an aperiodic and irreducible Markov chain that has the desired distribution as its equilibrium distribution \cite{Wea-PF}. The state of the chain after sufficient large number of simulation steps can be treated as an independent sample of the desired distribution. It is well known that MCMC method is an effective sampling method for complicated distributions in high-dimensional spaces. Taking the advantage of MCMC sampling, in this work we combine the solution of the backward SDE system with the observation data and apply the MCMC method to remove the stochastic space points away from statistically insignificant regions of the state pdf.

Specifically, in the Backward SDE filter framework, when we get the approximate solution $\Pi_{n}$ on $\mathcal{D}_{n}$ ($n =  1, 2, \dots, N_T-1$)  and initiate the recursive stage $t_n \rightarrow t_{n+1}$ by setting $P_n = \Pi_n$, instead of propagating $\mathcal{D}_n$ to $\mathcal{D}_{n+1}$ directly to construct the space point set for $\Pi_{n+1}$, we use MCMC sampling to create an observation informed intermediate point set $\mathcal{D}_{n+\f{1}{2}}$, and then propagate $\mathcal{D}_{n+\f{1}{2}}$ through the state model to get $\mathcal{D}_{n+1}$. To create the point set $\mathcal{D}_{n+\f{1}{2}}$, we generate a Markov chain for each space point $x_n^i \in \mathcal{D}_n$ to move it away from the statistically insignificant regions. Since the interpolatory approximation $\hat{P}_n$ in the Backward SDE filter is a point-wise numerical approximation of the state pdf, we use $\hat{P}_n$ as the stationary distribution of the Markov chain for each space point. It is worthy to point out that the Markov chain for each space point is based on the global approximation of the state pdf which is also incorporated with the observation information up to time level $t_n$. In this way, our MCMC resampling procedure is incorporated with the observation data, the data informed space points $\mathcal{D}_{n+\f{1}{2}}$ uses the observation information sufficiently and construct adaptive space points more effectively. There are many sampling algorithms for the MCMC method and we use Metropolis-Hastings algorithm\cite{Hastings} as an example to demonstrate our MCMC resampling method.

\subsubsection*{\underline{\textbf{Algorithm Summary.}}}

We summarize the recursive algorithms for Backward SDE filter for Jump Diffusion Processes as follows:

Define a pdf $Pr$ as the initial guess for the state $S_0$ by setting $\Pi_0 = Pr$, generate space points $\mathcal{D}_0 \sim \Pi_0$, and choose parameters $N$ as number of space points, $M$ as number of Monte Carlo samples, $J$ as number of meshfree interpolation points, $L$ as MCMC iteration number. The space point set $\mathcal{D}_1$ is propagated from $\mathcal{D}_0$ directly without resampling since there's no observation information at time $t = 0$.
\vspace{0.5em}

For the recursive stage $t_{n-1} \rightarrow t_{n}$, $n = 1, 2, \dots, N_T $, we implement
\vspace{-0.5em}
\begin{itemize}  \itemsep -2pt 
\item[-] \textit{Prediction Step:} Solve the backward SDEs system with scheme \eqref{Scheme} to get the predict state pdf $\tilde{P}_{n}$
\item[-] \textit{Update Step:} Update the state pdf with scheme \eqref{Bayes} to get $\Pi_{n}$ 
\item[-] \textit{Adaptive Meshfree Approximation:} Let $P_{n} = \Pi_{n}$ and expand $P_{n}$ on $\mathcal{D}_{n}$ to $\hat{P}_{n}$ through meshfree interpolation
\item[-] \textit{Resampling Step:} Use MCMC resampling procedure to construct observation informed space points $\mathcal{D}_{n+\f{1}{2}}$. Then propagate $\mathcal{D}_{n+\f{1}{2}}$ to $\mathcal{D}_{n+1}$ through the state model \eqref{NLF:State}
\end{itemize}

\subsection{Numerical examples}\label{Synthetic}

We present two numerical examples to demonstrate the performance of the Backward SDE filter in solving the nonlinear filtering problem \eqref{NLF}.  Both examples that we discuss in this subsection are benchmark problems with L\`evy noise added to the state equation as indicated in \eqref{NLF}.  In the first example, we estimate the trajectory of a target moves along a one dimensional periodic potential curve. In example 2, we solve a two dimensional bearing-only tracking problem \cite{zhang}.  We examine the performance of the Backward SDE filter by comparing our Backward SDE filter with auxiliary particle filter \cite{APF}, which is one of the most widely accepted nonlinear filtering methods by practitioners. The numerical experiments are carried out on an \textit{Intel Core i5 2.7 GHz} CPU.
 
\subsubsection*{Example 1}
\begin{figure}[h!]
\begin{center}
\includegraphics[scale = 0.4]{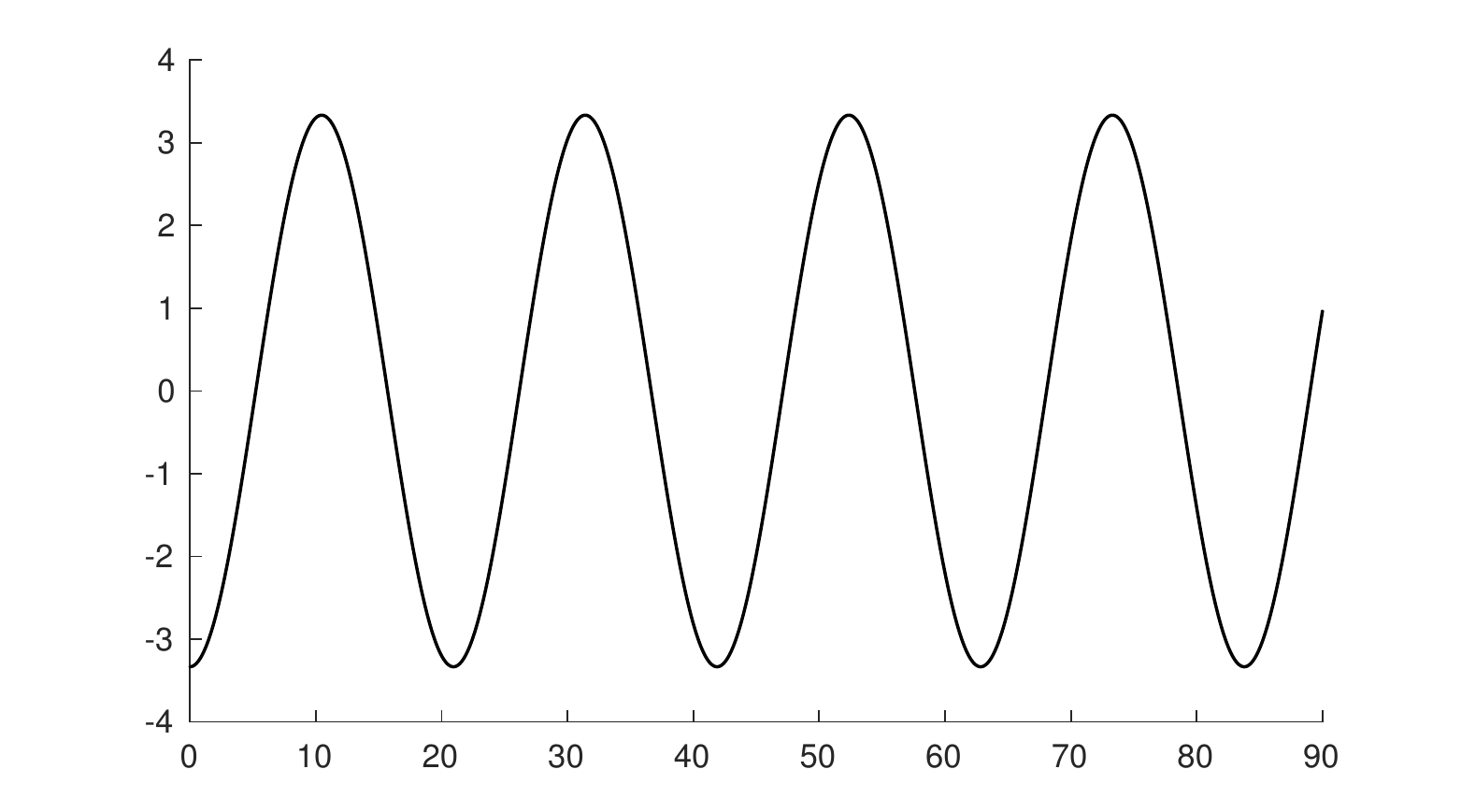} 
\end{center}
\caption{One dimensional periodic potential curve U}\label{1D_Potential}
\end{figure}

We first consider a periodic energy potential, denoted by $U$ with $\ds U = - \f{10}{3} \cos\Big(\f{3 x}{10}\Big)$, and assume that there's a target moves along the potential curve $U$. The curve is plotted in Fig. \ref{1D_Potential}.
We can see from this figure that the potential $U$ has some wells. If a particle moves on this potential curve $U$ and the potential difference forms the force to move the particle, the trajectory of the particle satisfies the ordinary differential equation $ \ds dS_t =  \sin\Big(\f{3 S_t}{10}\Big) dt$. Without the external perturbation, the particle wanders around the bottom of one of the wells. In this example, we assume that the particle is influenced by a L\`evy noise, which is the external perturbation to excite the particle from one potential well to another. The derived state process of the nonlinear filtering problem is given as follows,
\begin{equation*}  
dS_t = \sin\Big(\f{3 S_t}{10}\Big) dt + 4 dW_t +  \int_E 10 e \tilde{\mu} (dt, de ),
\end{equation*}
where $W_t$ is a standard Brownian motion and $\tilde{\mu}$ is the compensated Poisson measure. The dynamic system of this example is similar to the classic double well potential problem. But instead of two wells, it has multiple wells that allows a particle switch between.
The transitions of the particle and the evolution of the system is observed via data
given by the following observation process 
$$M_{t_n} = S_{t_n} + R \dot{B}_{t_n}, \quad n = 1 \cdots, N_t$$
where $B_t$ is a standard Brownian motion independent from $W_t$. In this problem, we track the particle as the target from time $t=0$ to $t =2$ with uniform time step $\Delta t = 0.02$, i.e. $N_t = 100$, and use a compound Poisson process to generate jumps in the state. For the observation, we choose $R = 0.1$. 
\begin{figure}[h!]
\begin{center}
\includegraphics[scale = 0.55]{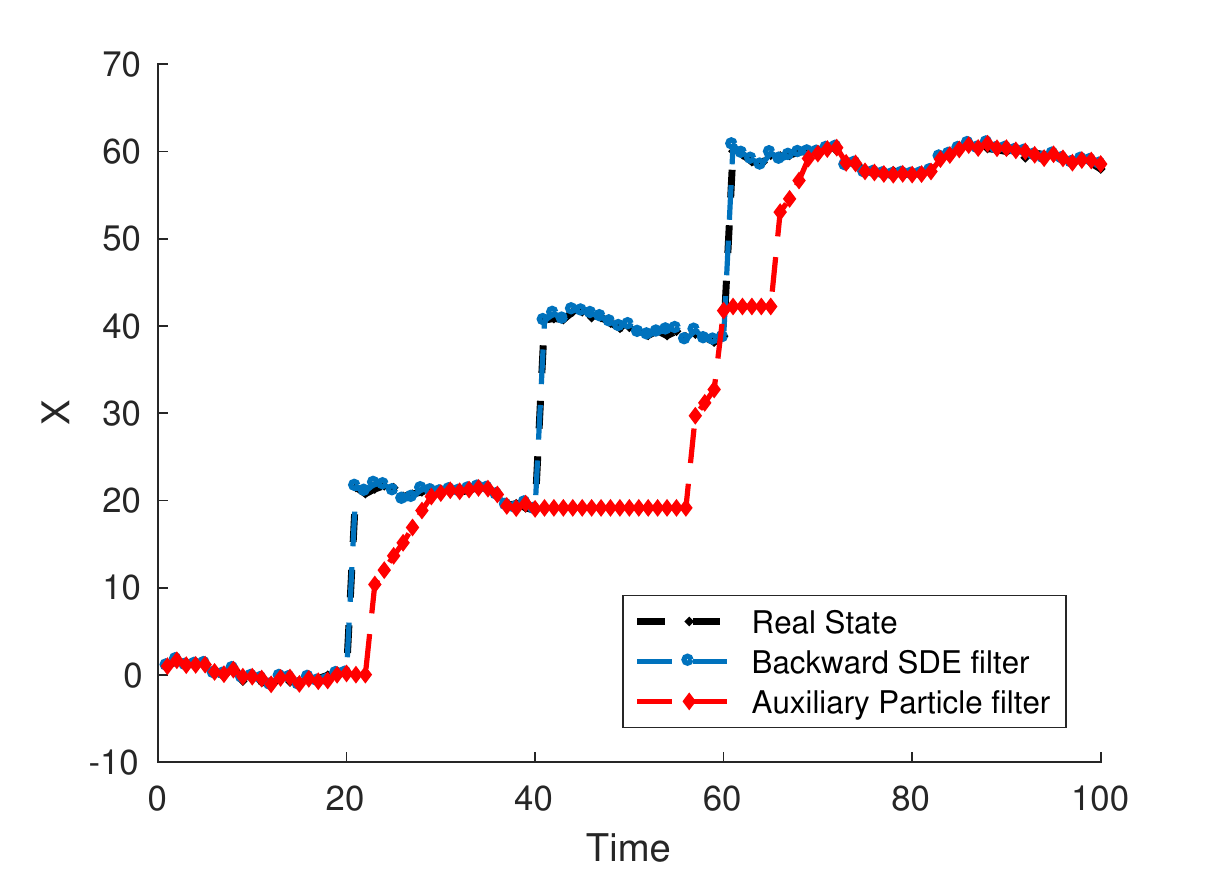} 
\end{center}
\caption{Comparison for one dimensional potential tracking. }\label{1D_Potential_Tracking}
\end{figure}
In Fig.  \ref{1D_Potential_Tracking}, we compare the performance of our Backward SDE filter with Auxiliary Particle filter (APF). In the Backward SDE filter, we use $200$ space points and in the APF, we use $800$ particles.
We can see from Fig. \ref{1D_Potential_Tracking} that the Backward SDE filter could capture the change of state very accurately while the auxiliary particle filter takes many more steps to capture the changes.
To further compare the performance of the Backward SDE filter and the APF, we repeat the above experiment $50$ times and compute overall root mean square errors (RMSE). In Table \ref{Efficiency}, we present the RMSE of the Backward SDE filter with $200$ space points as well as APF with $400$, $800$, $1600$ and $3200$ particles.
\renewcommand{\arraystretch}{1.25}
\begin{table} [h!]
\leftmargin=6pc \caption{Example 1: Efficiency comparison} \label{efficiency} \small
\begin{center}\label{Efficiency} 
\begin{tabular}{|c|c|c|}
 \hline  Numerical methods &  CPU time (seconds) & $err_{G}$   \\
\hline   Backward SDE filter ($200$ space points) &   $9.01$  & $0.3302$\\
\hline   Auxiliary particle filter ($400$ particles) &  $10.97$ &  $34.1392$\\
\hline   Auxiliary particle filter  ($800$ particles) &  $ 23.07$ & $27.5474$\\
\hline   Auxiliary particle filter  ($1600$ particles) &  $ 50.49$  & $16.7304$\\
\hline   Auxiliary particle filter  ($3200$ particles) &  $ 95.40$  & $10.8136$\\
\hline
\end{tabular}\end{center}
\end{table}
From this Table, we can see that with $200$ space points, the Backward SDE filter takes $9$ seconds to tack the target trajectory and the RMSE is $0.3302$. Although the APF with $400$ particle could track the target with similar CPU time to the Backward SDE filter, the RMSE of APF is $34.1392$, which is much higher than RMSE of the Backward SDE filter. When using more and more particles in APF, the RMSE reduces and the computing cost increases significantly. However, even using $3200$ particles to track the target in the APF, which consumes more than $10$ times of CPU time of the Backward SDE filter, the RMSE of APF is still much higher than the Backward SDE filter.

\subsubsection*{Example 2}
In this example, we solve a bearing-only tracking problem, in which a target is moving one a two dimensional plane with a near constant velocity and the state of the target is perturbed by L\`evy noise. The state equation is given as follows 
\begin{equation}\label{Ex2:state}
d \bm{S}_{t} = \bm{A} \bm{S}_{t} dt + \bm{\sigma} dW_{t} +  \int_E  \bm{\beta(e)} \tilde{\mu} (dt, de ),
\end{equation}
where $\bm{S}_{t} = (X_{t}, Y_{t}, \dot{X}_{t}, \dot{Y}_{t})^{T}$ is a 4 dimensional vector, $(X_{t},  \dot{X}_{t})$ and $(Y_{t},  \dot{Y}_{t})$ are the position and velocity of the target  at time $t$ corresponding to $X$ and $Y$ co-ordinate, respectively. $W_{t}$ is a 4 dimensional standard Brownian motion and $\tilde{\mu}$ is the compensated Poisson measure. The matrices $A$, $\bm{\sigma}$ and $\bm{\beta}$ are given by
$
\bm{A} = \left(
  \begin{array}{cc}
    I_2 & I_2 \\
    0 & I_2
      \end{array}
\right),  \hspace{1em}
\bm{\sigma} = diag(0.1, 0.1, 0.05, 0.05), \hspace{1em} 
\bm{\beta(e)} = \bm{e} (2, 2, 0.2, 0.2)^{T},
$
where $I_2$ is $2 \times 2$ identity matrix.
The target is observed by a bearing-range director located at $
( X_{obs}, Y_{obs})$, i.e. the measurement $\bm{M}_t$ is given by
\begin{equation*}\label{Ex2:measurement}
\bm{Y}_{t} = \left(
\begin{array}{ccc}
 \arctan\left( \ds \f{Y_{ t} - Y_{obs}}{ X_{t} - X_{obs} } \right) \\
 \sqrt{( X_{ t_n} - X_{obs})^2 + (Y_{t} - Y_{obs})^2}
\end{array}\right) + \bm{R} \dot{B}_{t}, \quad n = 1, 2, \cdots, N_t,
\end{equation*}
where $B_{t}$ is a standard two dimensional Brownian motion independent from $W_{t}$ and $\bm{R} := diag(0.01, 0.1)$ is a two dimensional matrix.

In this example, we use $\alpha$-stable process to generate jumps in the noise and track the target for the time period $t = 0$ to $t =2$ with uniform time step $\Delta t = 0.04$, i.e. $N_t = 50$.  
\begin{figure}[h!]
\begin{center}
\includegraphics[scale = 0.55]{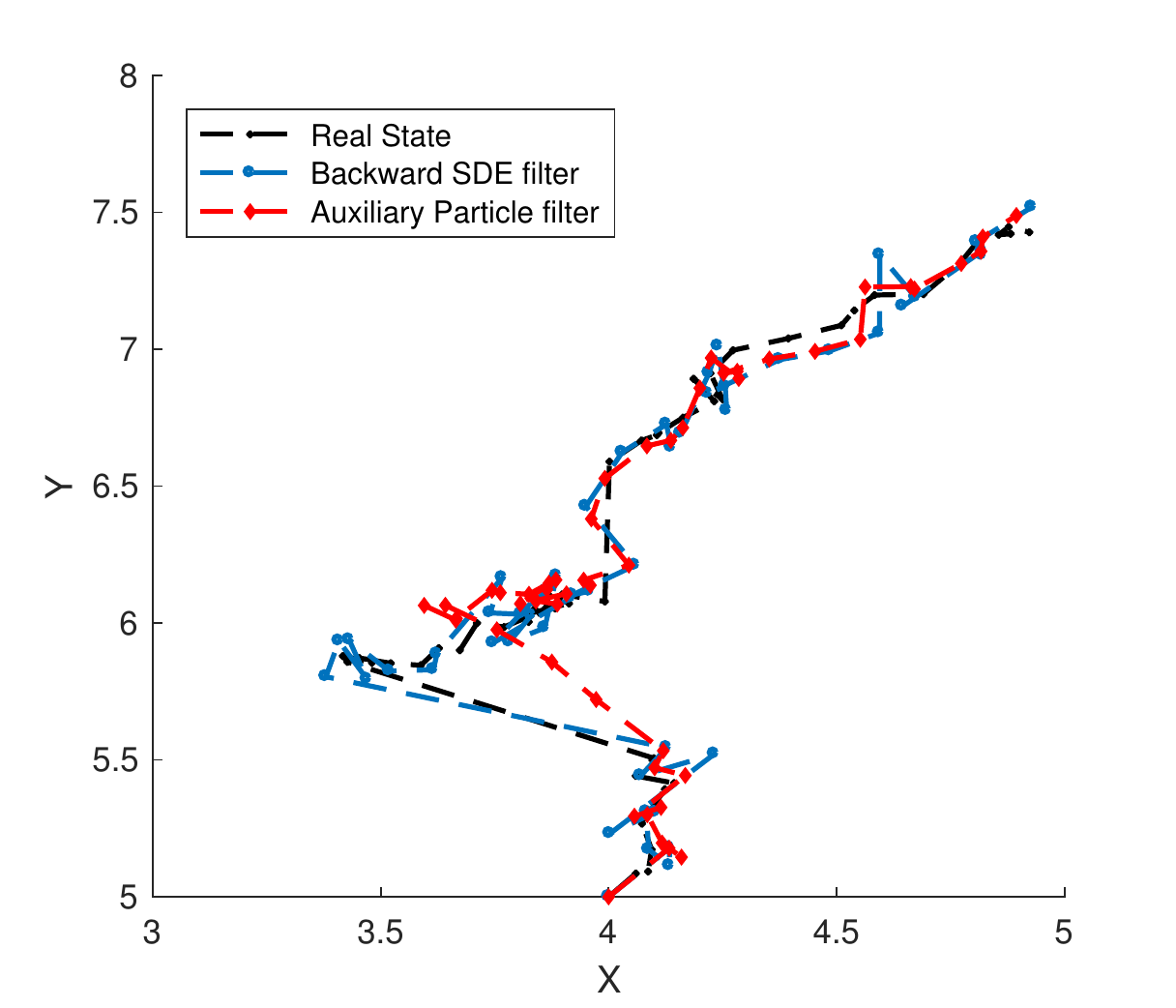} 
\end{center}
\caption{Comparison of target tracking for $\alpha = 1$ }\label{4D_Tracking}
\end{figure}
In Fig. \ref{4D_Tracking}, we compare the tracking performance of the Backward SDE filter with APF and choose $\alpha =1$. For the Backward SDE filter, we use $1,500$ space points. For the APF, we use $6,000$ particles to describe the conditional pdf of the target state. The black curve is the real target trajectory in the $XY$-plane. We can see there are several jumps in the trajectory. The red curve in the figure is the estimate trajectory obtained by the APF, and the blue curve is the estimate trajectory obtained by the Backward SDE filter.
We can see from Fig. \ref{4D_Tracking} that the Backward SDE filter could capture the change of state much faster than the auxiliary particle filter.

To better demonstrate the effectiveness of our Backward SDE filter in tracking a target with jumps, we choose $\alpha = 0.5$ for the $\alpha$-state L\`evy process and solve the same bearing-only tracking problem. \begin{figure}[h!]
\begin{center}
\includegraphics[scale = 0.55]{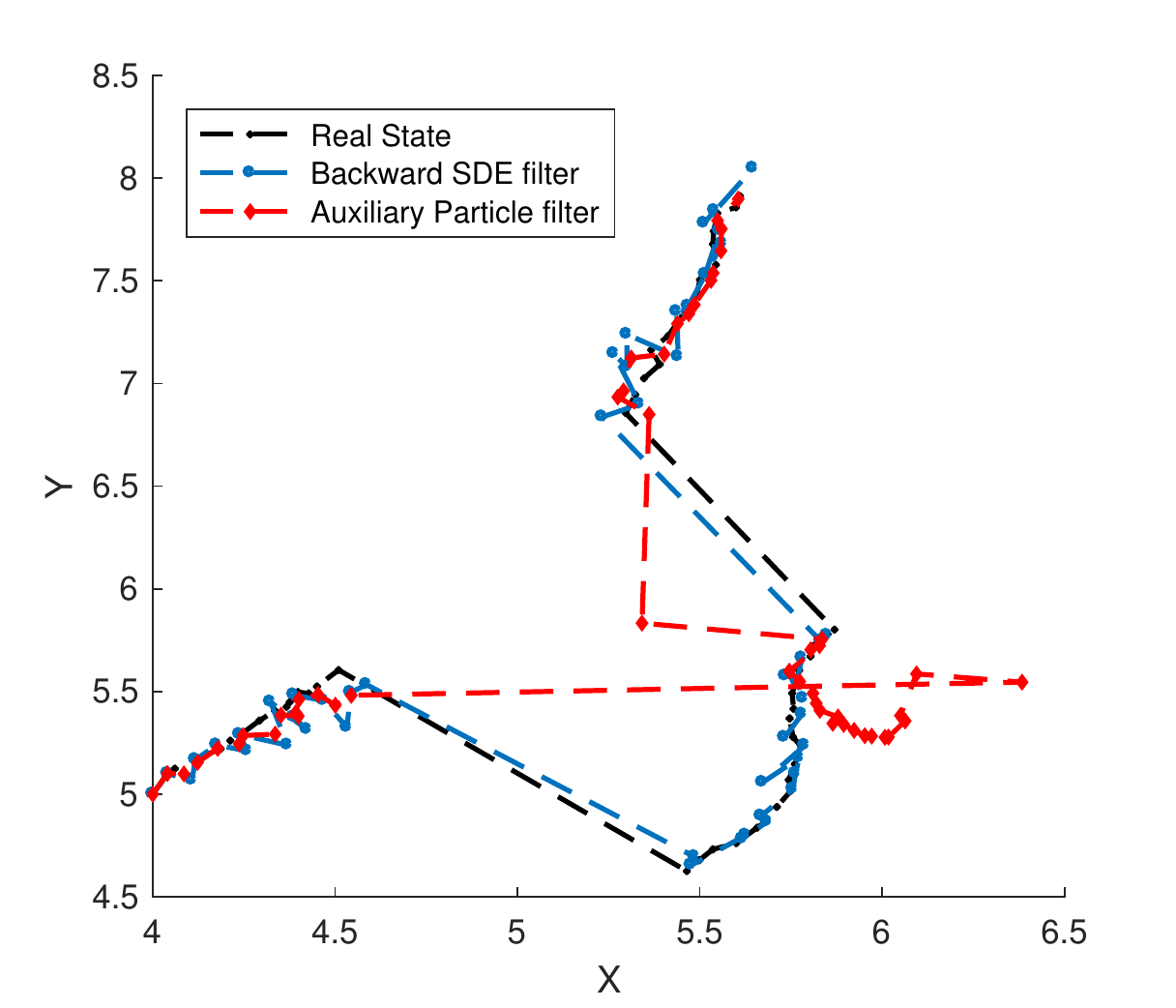} 
\end{center}
\caption{Comparison of target tracking for $\alpha = 0.5$ }\label{4D_Tracking_2}
\end{figure}
Since we choose a small $\alpha$ in this experiment, the $\alpha$-stable process has thicker tail in the distribution, which results higher probability for big jumps. 
In Fig. \ref{4D_Tracking_2}, we present the tracking performance for the case $\alpha = 0.5$ and we can observe two significant jumps in the target trajectory. 
From this figure, when choosing a smaller $\alpha$, we can see that estimate of APF lost target for several steps and the irregular behavior when trying to catch up with the target is typically caused by its degeneracy problem. However, our Backward SDE filter still provides reliable estimate for the target state even with two big jumps.

\section{Application of Backward SDE filter in material science }\label{Material}

In this section, we discuss the application of the Backward SDE filter in nano-phase material sciences, in the context of development of algorithmic control over single atom processes that will enable large-scale and automated manipulation and synthesis of single atoms and molecules  \cite{Hla_2014}, \cite{Kim_2015}, \cite{Morgenstern_2013}. The majority of studies on atomic-scale manipulation to date have been carried out largely in the open-loop regime, where there is no direct feedback between excitation (electric field, electric current, direct chemical interaction) and the excited object, or in rare but intriguing cases -- with human control \cite{Green_2014}. However, automated control methods will be needed to achieve two major goals: dramatically increase the manipulation speed, well beyond the human capacity, and equally importantly to enable deterministic selectivity over the reaction steps -- a specific kind of reaction, specific direction of motion etc \cite{Peter_book}.  

The automated synthesis algorithm can be generally described by three procedures: the stochastic optimization procedure, which designs the optimal material potential surface; the tracking procedure, which tracks the movement of a target atom based on the designed material potential surface and the observation data received from scanning probe (or related electron microscopy); the optimal control procedure, which controls the material condition to minimize the cost of synthesis based on the estimation of the molecule state. In this atomic level material synthesis, the role of the Backward SDE filter is to track the atom trajectory in the single-atom manipulation and single-molecule reaction, both of which are registered as abrupt events in the relevant observables, such as tunneling current and/or interaction force between the manipulating probe and the manipulated entity.

To examine the applicability of our algorithm in the aforementioned framework, we use Backward SDE filter to simulate tracking atoms moving on the potential surface of some well-known materials. Single atoms on the surface experience atomic-scale interaction potential, due to preferential atomic bonding toward maximally coordinated sites  \cite{Morgenstern_2013}. In this work, we approximate the corresponding potential using a simple triangular lattice of potential energy wells as a sum of sinusoids, where the amplitude was calibrated to match the experimentally observed diffusion potential. This would coarsely resemble the potential of the 111-terminated noble-metal surface (such as Ag(111), Au(111) and others) toward interaction with a single atom. The depth of the wells was chosen to mimic a recent work on atomic motion by Giessibl et al.  \cite{Ternes_2008}. In what follows, we present the performance of the Backward SDE filter in tracking the atom trajectories on approximated potential surfaces.

\vspace{0.25em}

\textbf{Experiment 1.}
\vspace{0.25em}

\begin{figure}[h!]
\begin{center}
\subfloat[3D image of the material surface potential. ]{\includegraphics[scale = 0.4]{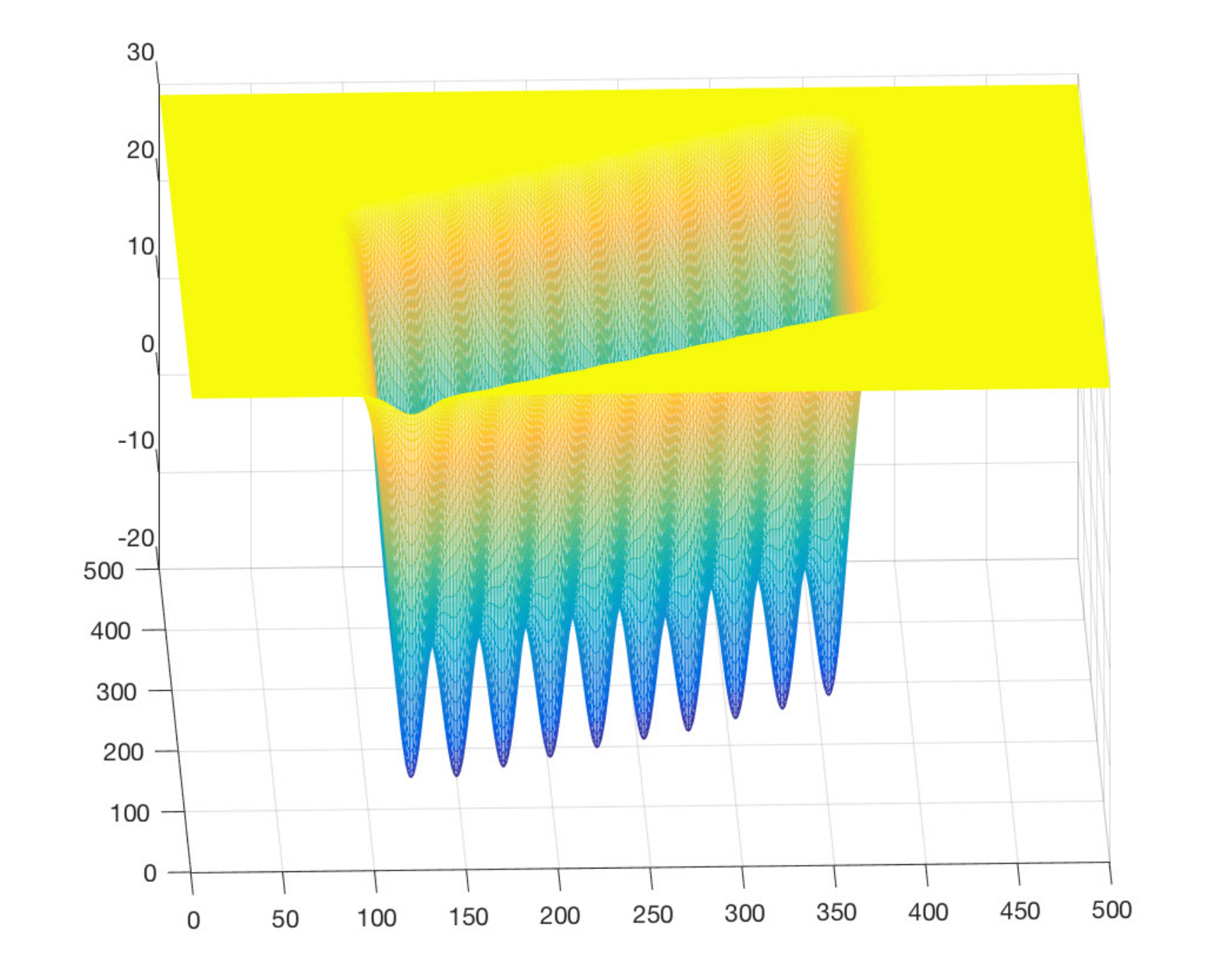} }
\subfloat[2D plan view of the material surface potential. ]{\includegraphics[scale = 0.38]{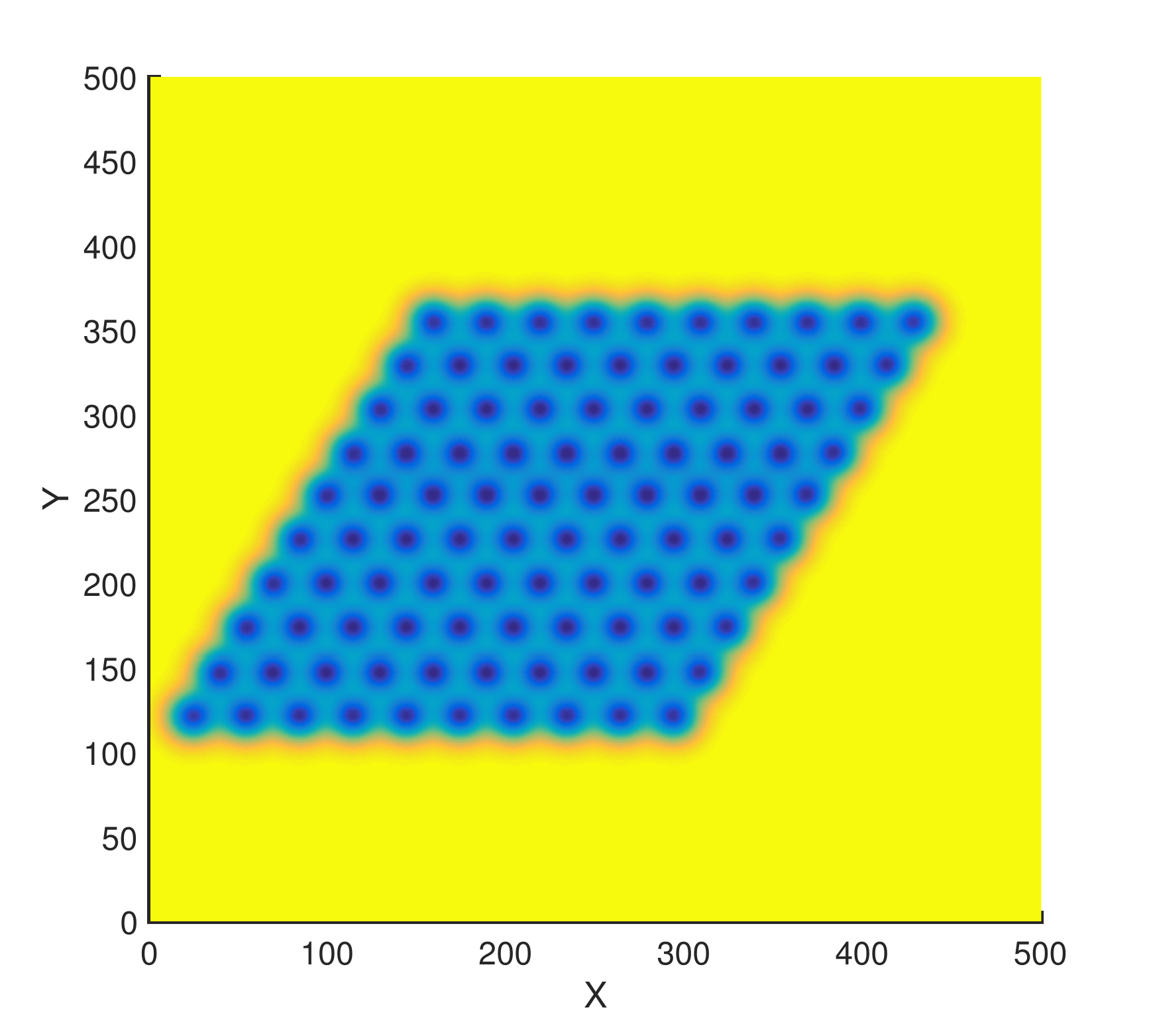} }
\end{center}
\caption{The material surface potential $F_1$.}\label{Real_Potential}
\end{figure}

We first depict our approximated potential in Fig. \ref{Real_Potential} (a), which represents the potential energy landscape for diffusion of atomic scale species (such as atoms or molecules), and denote this potential by $F_1$. From this image, we can see that there are several deep wells. Once the atom falls in the bottom of one well, it will be trapped in this well unless some excitation occurs.
In Fig. \ref{Real_Potential} (b), we present the 2D plan view of the energy potential $F_1$. The dark blue disks represent low potential regions which are bottoms of wells and the connected light blue region represent high potential area.

\begin{figure}[h!]
\begin{center}
\subfloat[Trajectory in X coordinate]{\includegraphics[scale = 0.5]{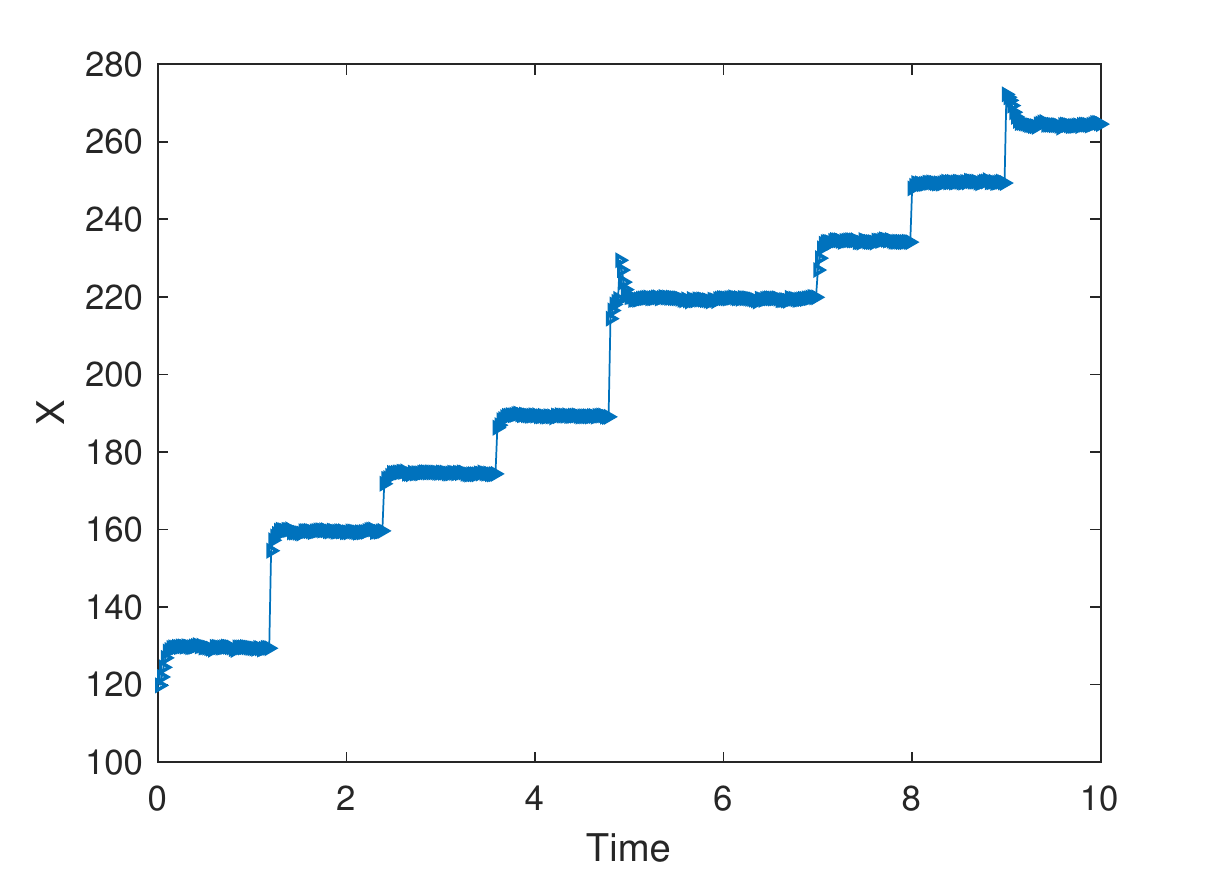} }
\subfloat[Trajectory in Y coordinate]{\includegraphics[scale = 0.5]{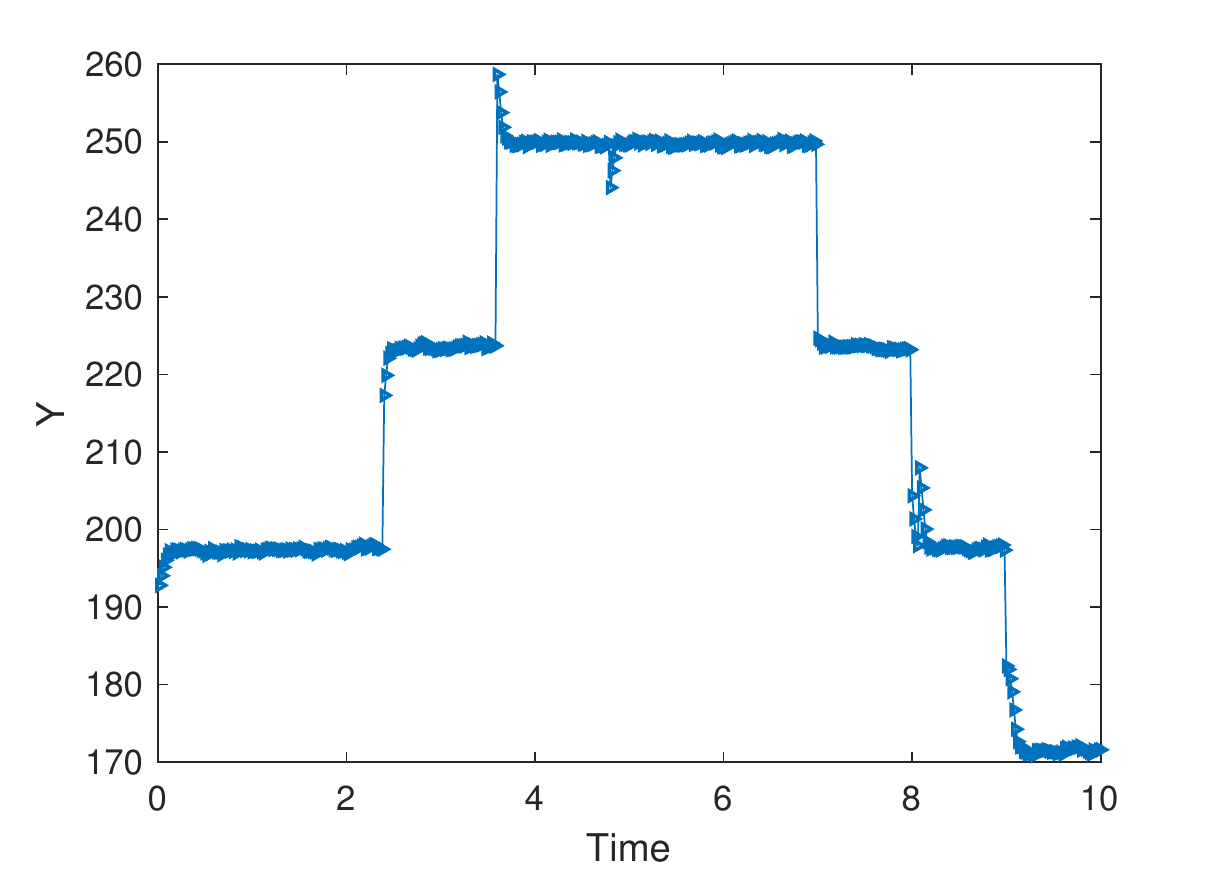} }
\end{center}
\caption{ Sample trajectory in X, Y directions} \label{XY_Potential}
\end{figure}
\begin{figure}[h!]
\begin{center}
\includegraphics[scale = 0.48]{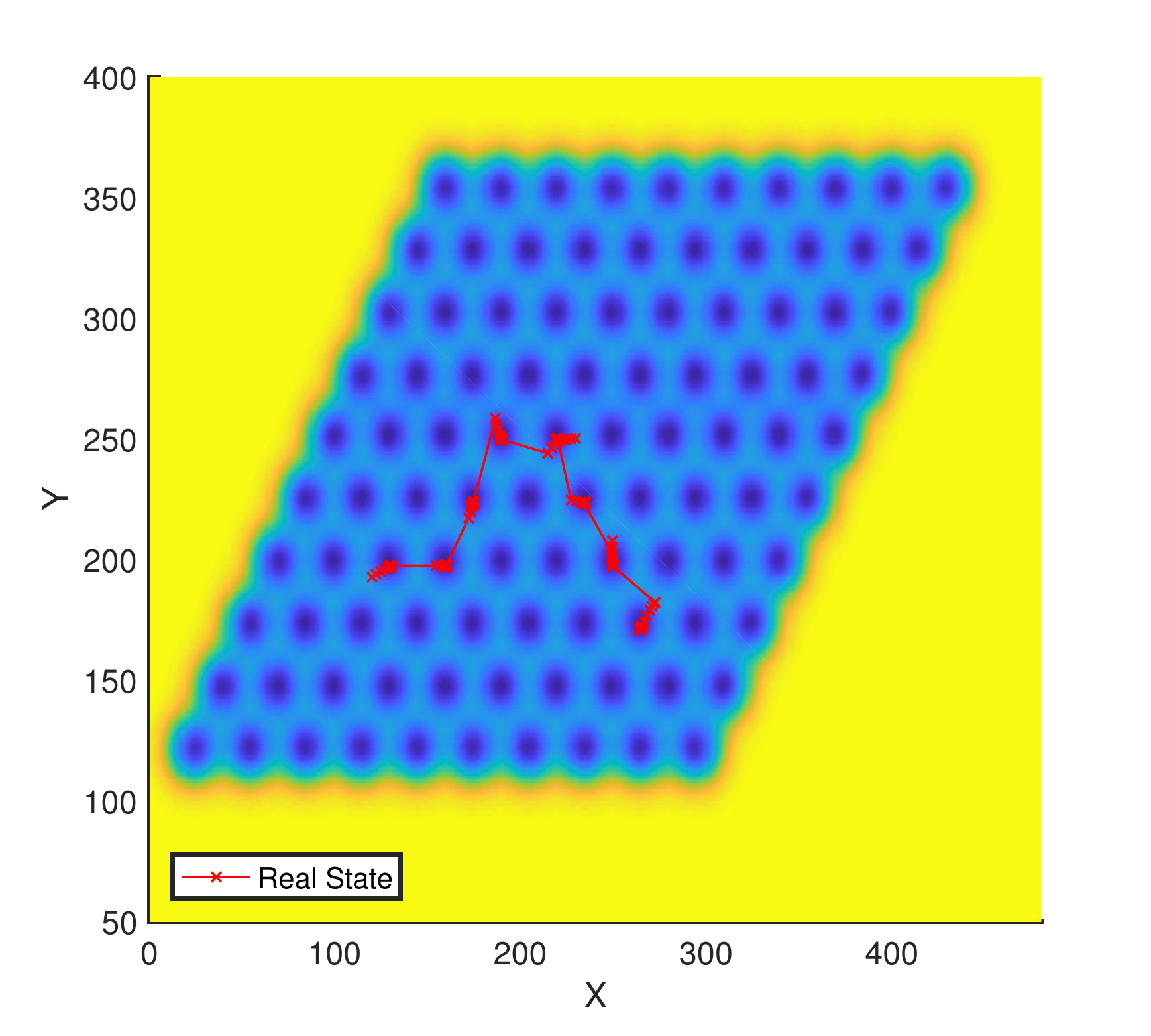} 
\end{center}
\caption{The sample atom trajectory with respect to the energy potential $F_1$ }\label{Potential_Trajectory}
\end{figure}

The position of an atom would follow the force caused by the potential of material surface presented in Fig. \ref{Real_Potential}. Different from the synthetic examples we presented in Section \ref{Synthetic}, in which we have analytic drift function $b$ in the nonlinear filtering problem \eqref{NLF}, the drift term of the state equation in this experiment is calculated by the simulated physical potential force $F_1$.
\begin{figure}[h!]
\begin{center}
\includegraphics[scale = 0.55]{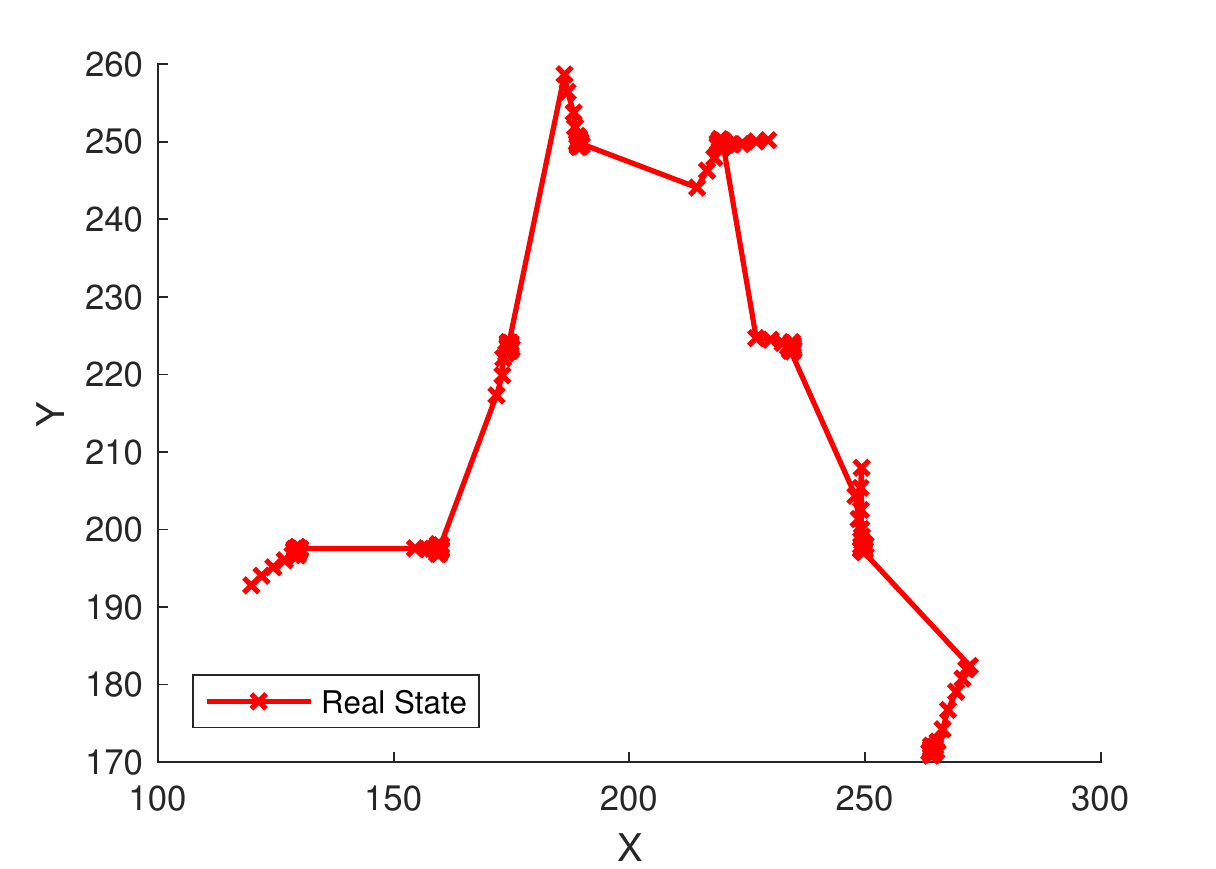} 
\end{center}
\caption{A sample atomic trajectory}\label{Potential_XY}
\end{figure}
As a result, the state in this experiment is given by
\begin{equation}\label{state:potential1}
dS_t = - \nabla \hat{F_1}(S_t) dt + \bm{\sigma} dW_t + \int_E \bm{\beta(e)} \tilde{\mu}(dt, de),
\end{equation}
where $S_t = (X_t, Y_t)^T$ is a two dimensional vector describing the position of the atom in the $XY$ plane. In this experiment, the function $\nabla \hat{F_1}$ is the drift term of the state equation given by the force caused by the potential surface. Since the energy potential $F_1$ is obtained from simulation on a given grid mesh, it is different from the synthetic examples in Section \ref{Synthetic}, in which the drift term $b$ is an explicit function. In order to derive the gradient of $F_1$, we use polynomial approximation to construct a smooth surface of $F_1$ and then calculate the gradient of the surface, i.e. $\nabla \hat{F_1}$ on the $XY$ plane. To simulate the trajectory, we choose $\bm{\sigma} = diag(0.1, 0.1)$, $\bm{\beta(e)} = \bm{e} (10, 10)^{T}  $ and use a compound Poisson process to generate jumps in the state.
In Fig. \ref{Potential_XY}, we present a sample atom trajectory of the simulated state $S_t$ over time interval $[0, 10]$ with uniform time step $\Delta t = 0.02$, i.e. $N_t = 500$.
From this figure, we can see that the atom has some jumps with several stable states. Once the atom arrives at a stable state, it remains at this state until excited by some force and jumps to another random state. To demonstrate the details of the atom trajectory, we plot the atom position in each direction in Fig. \ref{XY_Potential}. From this figure, we can see that there are $8$ stable positions, in which the atom remains in the similar $X-Y$ position for a period, and the atom has some random smaller scale jumps which is not enough to escape from the bottom of the potential well and is dragged back by the potential force.
In Fig. \ref{Potential_Trajectory}, we plot the sample atom trajectory depicted in Fig. \ref{Potential_XY} in the energy potential which has been shown in Fig. \ref{Real_Potential} (b).
From this figure, we can see that all the stable positions of the atom trajectory is around the bottom of a potential well.

In order to track this atom, we use tunnel electron microscope to receive observation of the atom. In this experiment, we assume that the observation is noise perturbed atom position, i.e. 
$$M_{t_n} = S_{t_n} + R \dot{B}_{t_n}, \quad n = 1,2,\cdots, 500,$$ 
\begin{figure}[h!]
\begin{center}
\includegraphics[scale = 0.6]{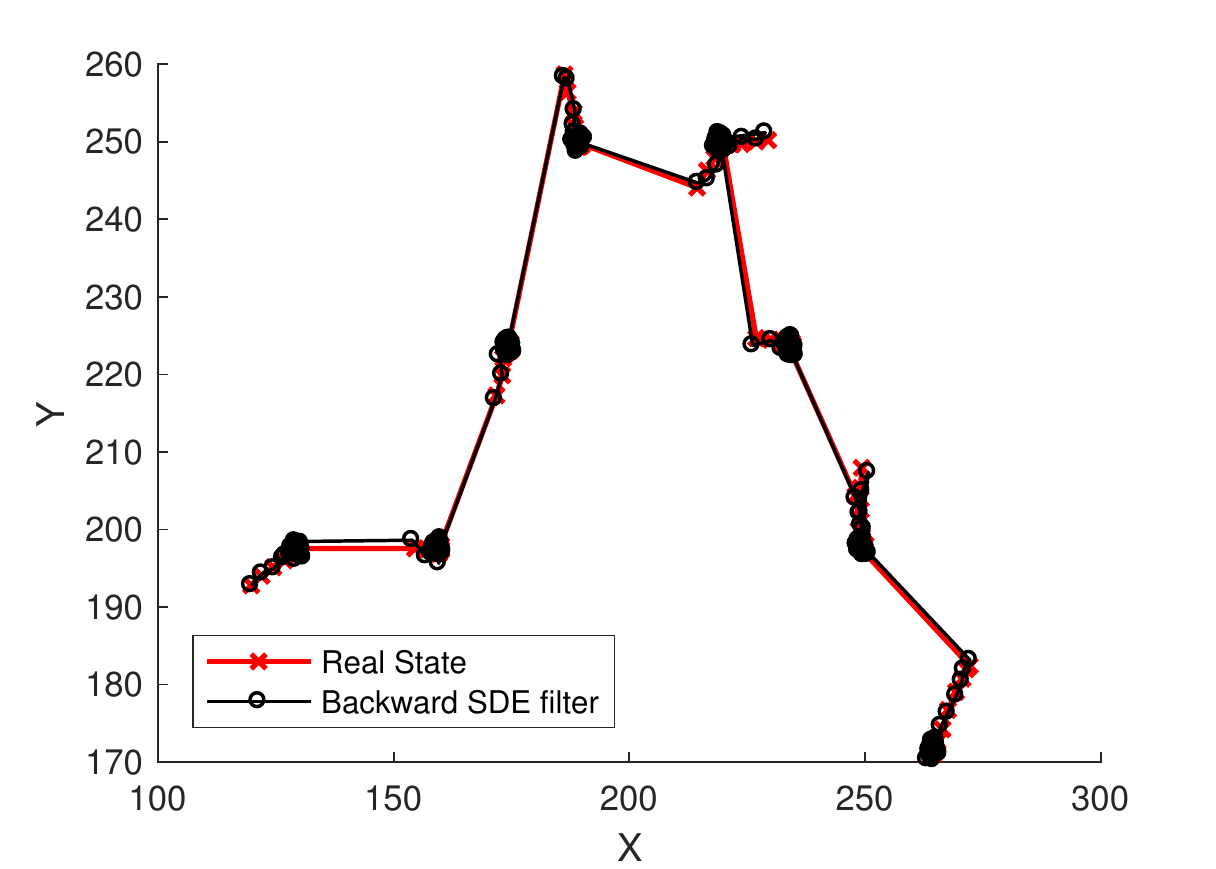} 
\end{center}
\caption{Tracking performance of Backward SDE filter in the real energy potential problem}\label{Tracking_Potential}
\end{figure}
where $B_t$ is a two dimensional Brownian motion independent from $W_t$ and we choose $R = diag(0.05, 0.05)$.
In Fig. \ref{Tracking_Potential}, we present the performance of our Backward SDE filter in tracking this atom . 
The red trajectory in this figure is the real target atom trajectory and the black path marked by circles is the estimated target obtained by using the Backward SDE filter. We can see from this figure that the Backward SDE filter could track the atom trajectory accurately.
\vspace{0.25em}

\textbf{Experiment 2.}
\vspace{0.25em}

In this experiment, we approximate a different energy potential surface $F_2$ and use the Backward SDE filter to track the atom trajectory based on observation $M_{t_n}$. 
\begin{figure}[h!]
\begin{center}
\subfloat[3D image of the material surface potential.]{\includegraphics[scale = 0.4]{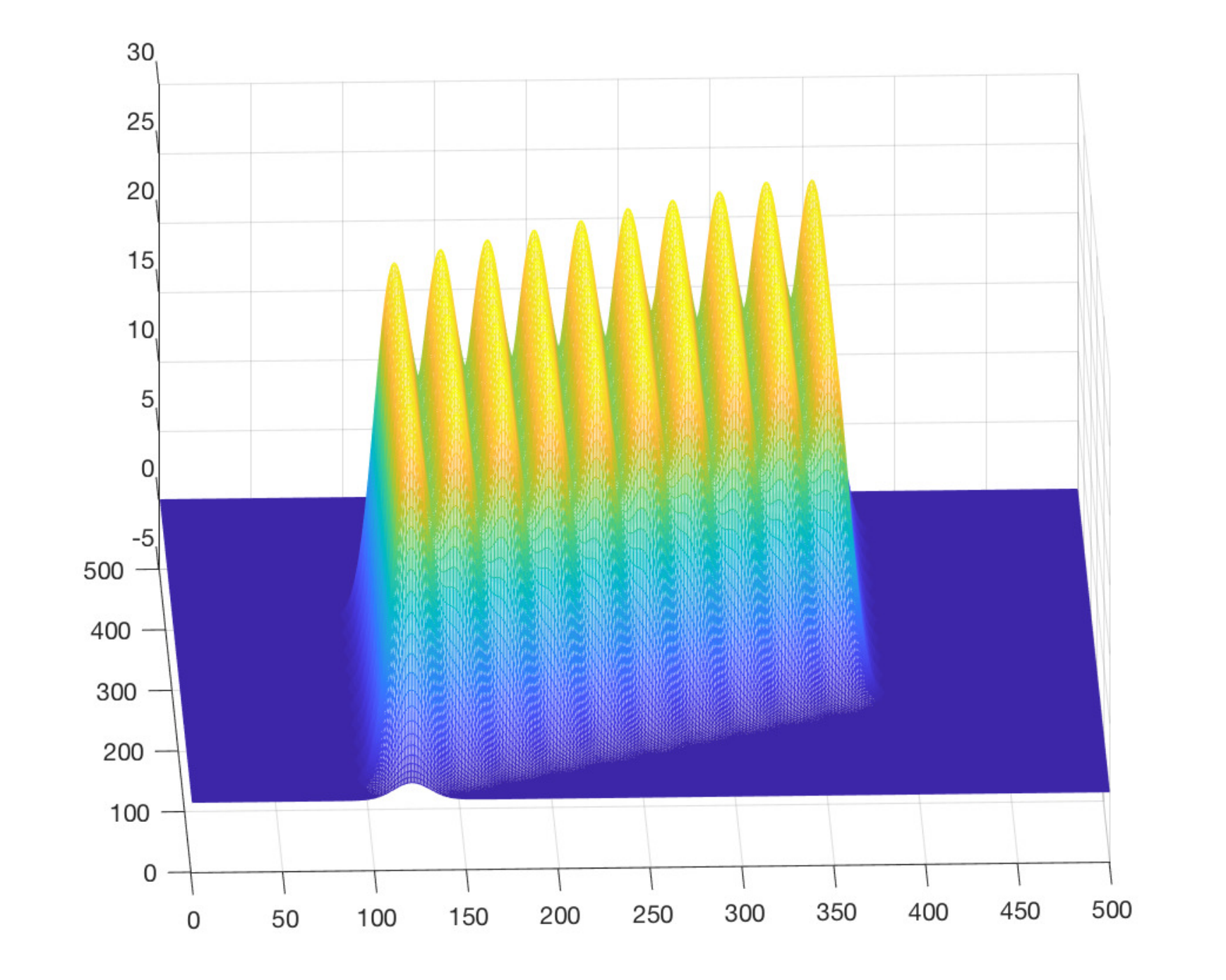} }
\subfloat[2D plan view of the material surface potential.]{\includegraphics[scale = 0.38]{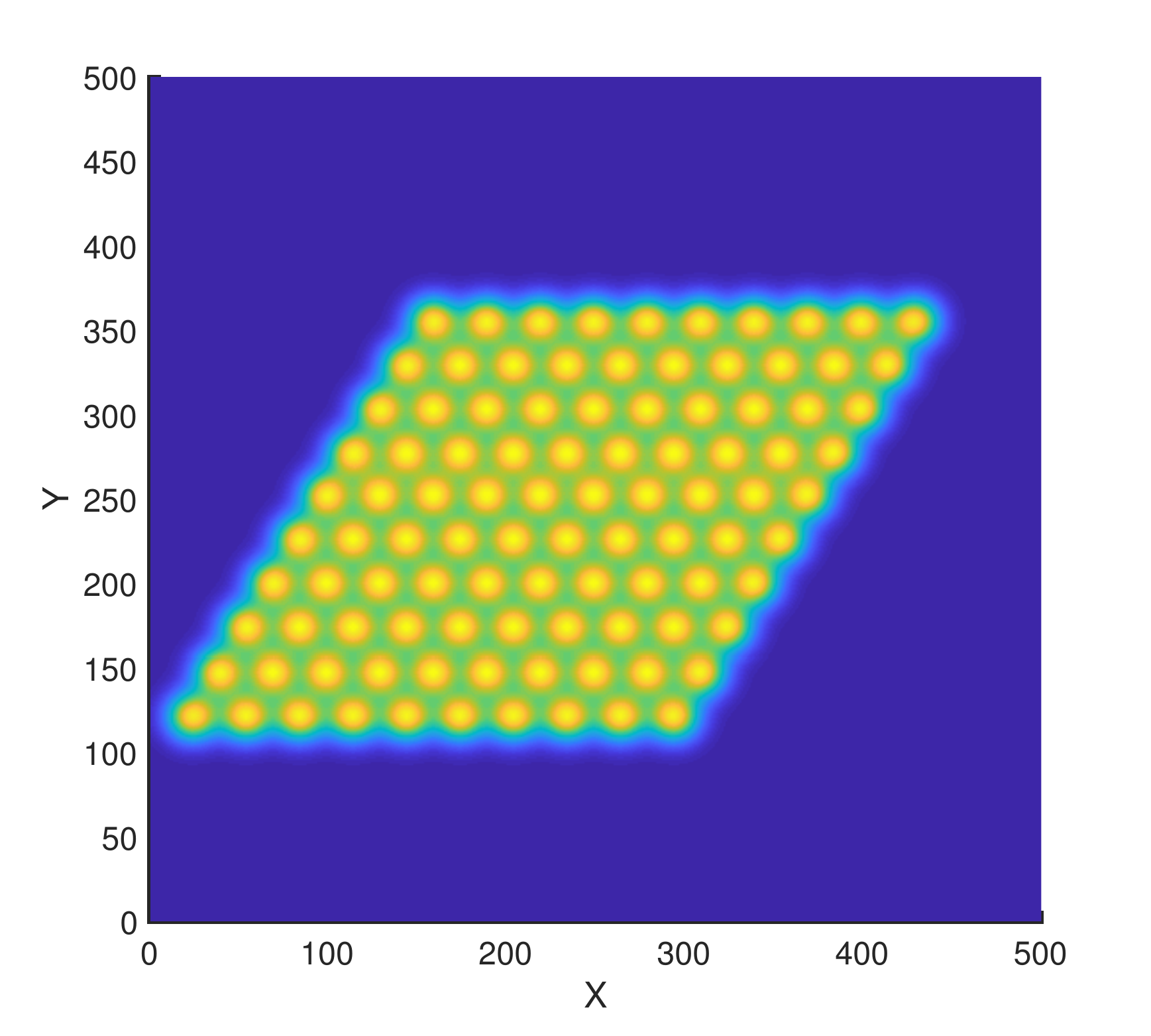} }
\end{center}
\caption{The material surface potential $F_2$. }\label{Real_Potential2}
\end{figure}
In Fig. \ref{Real_Potential2} (a), we present the 3D potential surface that we approximated. On the contrary to energy $F1$, we can see there are several peaks for the energy potential. If the atom gets on the peak of potential by any excitation, it would slide down and remains in the bottoms of different peaks which has lower energy with more stable state.
In Fig. \ref{Real_Potential2} (b), we plot the 2D plan view image of the energy potential $F_2$. The yellow disks represent the energy peaks and the green region has lower energy in which the atom is in stable state. Given the simulated energy potential $F_2$, we derive the state equation of the nonlinear filtering problem \eqref{NLF} as
\begin{equation*}
dS_t = - \nabla \hat{F_2}(S_t) dt + \bm{\sigma} dW_t + \int_E \bm{\beta(e)} \tilde{\mu}(dt, de),
\end{equation*}
where $\nabla \hat{F_2}$ is the drift term of the state equation, which is the approximate gradient of the energy potential surface $F_2$.

\begin{figure}[h!]
\begin{center}
\includegraphics[scale = 0.55]{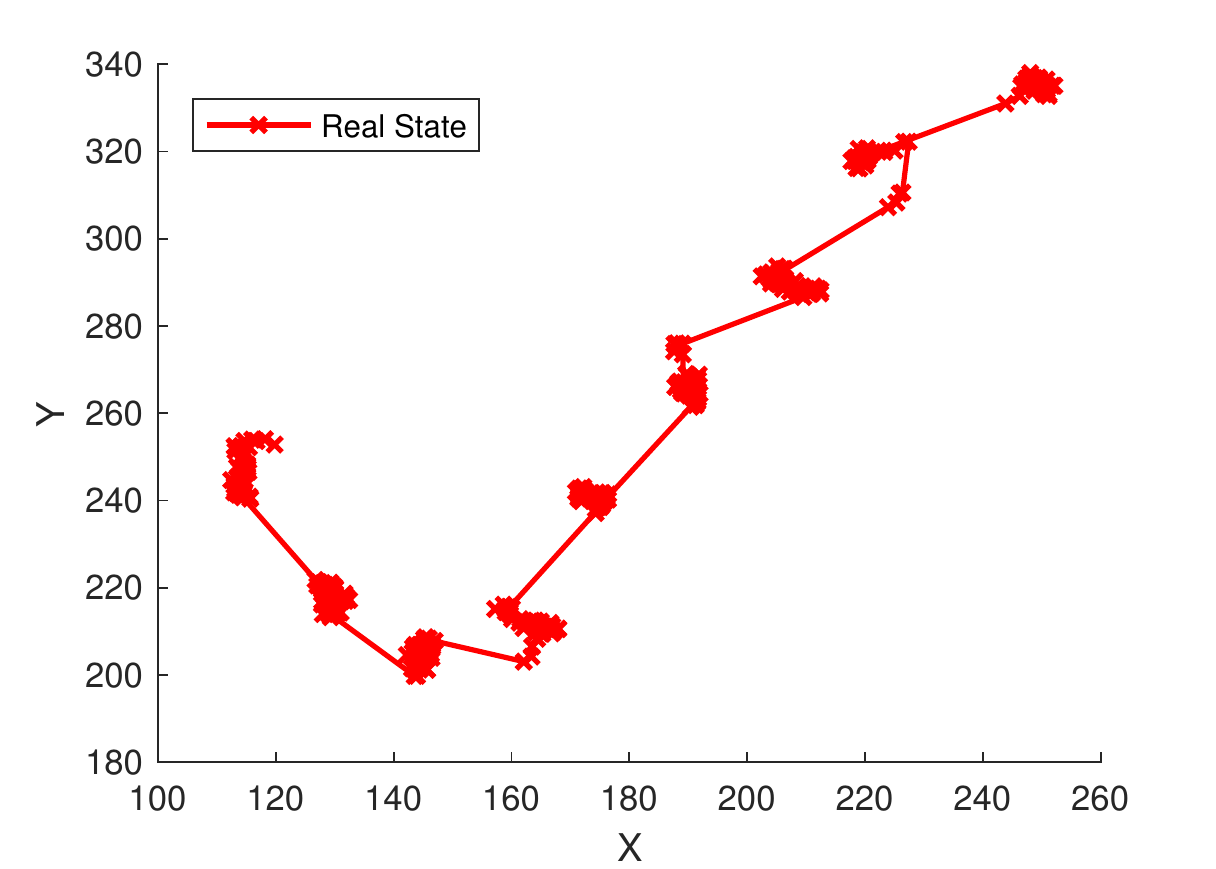} 
\end{center}
\caption{A sample atomic trajectory}\label{Potential2_XY}
\end{figure}
\begin{figure}[h!]
\begin{center}
\includegraphics[scale = 0.48]{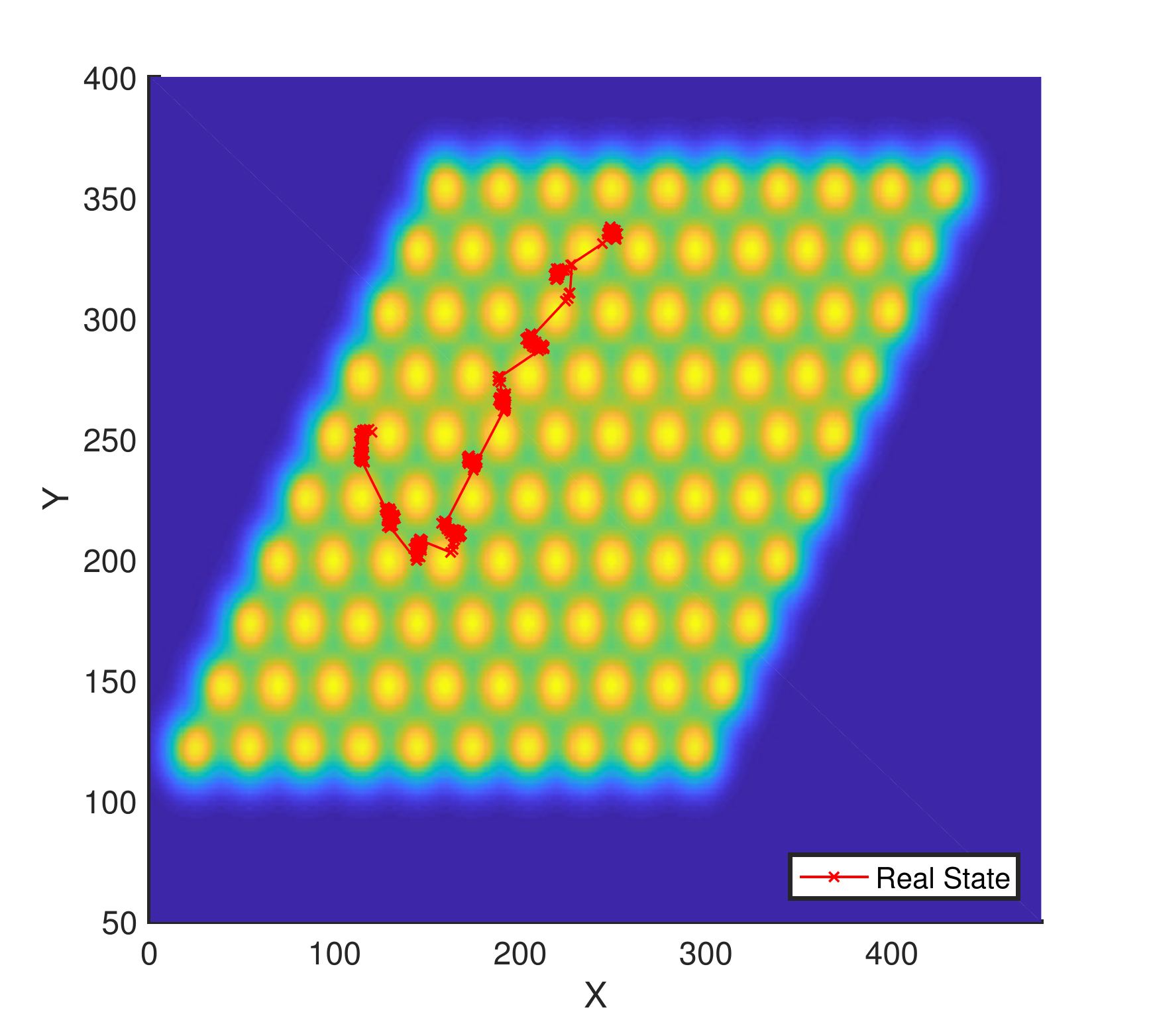} 
\end{center}
\caption{The sample atomic trajectory with respect to the energy potential $F_2$.}\label{Potential2_Trajectory}
\end{figure}

Similar to Experiment 1, we simulate the state equation over time interval $[0, 10]$ with uniform time step $\Delta t = 0.02$, i.e. $N_t = 500$ with $\bm{\sigma} = diag(0.1, 0.1)$ and $\bm{\beta(e)} = \bm{e} (10, 10)^{T}  $. In Fig. \ref{Potential2_XY}, we plot a sample atom trajectory and put this trajectory on the material surface in the $XY$ plane with the 2D plan view of energy potential $F_2$ in Fig. \ref{Potential2_Trajectory}. 
From Fig. \ref{Potential2_Trajectory}, we can see that for most of time, the atom remains stable in the lower energy region.  The position of the atom is perturbed by both Gaussian noises, which makes it linger around its current position, and Poisson noise, which causes some jumps. If the atom jumps onto one of the peak, as we can see in the figure, the potential force pushes it down and keep it remain in the bottom of all the peaks.

\begin{figure}[h!]
\begin{center}
\includegraphics[scale = 0.6]{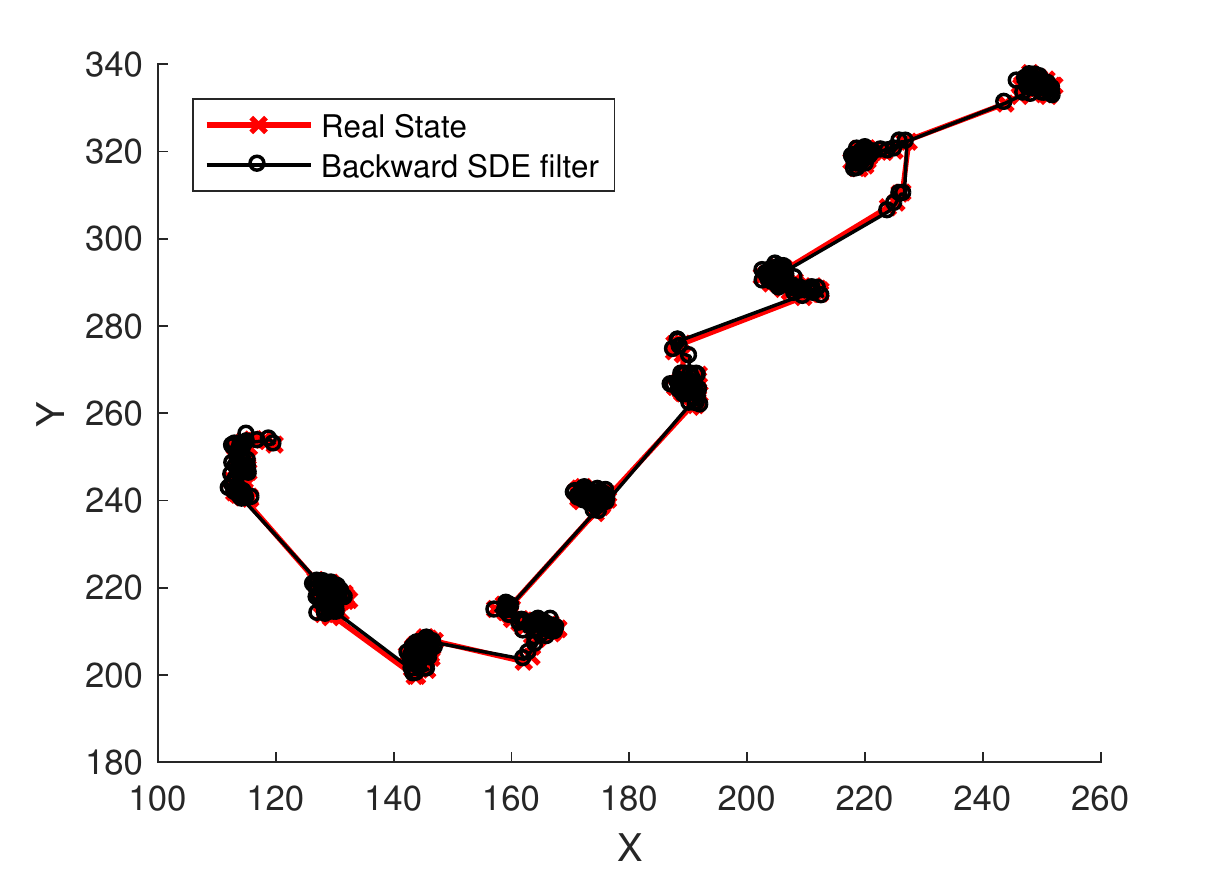} 
\end{center}
\caption{Tracking performance of Backward SDE filter in the real energy potential problem}\label{Tracking_Potential2}
\end{figure}
Finally, we present the performance of our Backward SDE filter in tracking this atom in Fig. \ref{Tracking_Potential2}. 
The red trajectory in this figure is also the real target atom trajectory and the black path marked by circles is the estimated target obtained by using the Backward SDE filter. We can see from this figure that under the energy potential $F_2$, the Backward SDE filter could also track the atom trajectory accurately.

\section{Acknowledgement}

This work is partially supported by U.S. Department of Energy's Advanced Scientific Computing Research through the ACUMEN project. The second author acknowledges support by U.S. National Science Foundation under grant number DMS-1720222 and the third author also acknowledges support by the Center for Nanophase Materials Sciences, sponsored by the Division of User Facilities, Basic Energy Sciences, U.S. Department of Energy.




\begin{thebibliography}{10}

\bibitem{MCMC-PF}
C.~Andrieu, A.~Doucet, and R.~Holenstein.
\newblock Particle {M}arkov chain {M}onte {C}arlo methods.
\newblock {\em J. R. Stat. Soc. Ser. B Stat. Methodol.}, 72(3):269--342, 2010.

\bibitem{BaoC20142}
F.~Bao, Y.~Cao, and X.~Han.
\newblock Forward backward doubly stochastic differential equations and optimal
  filtering of diffusion processes, arxiv: 1509.06352.
\newblock 2016.

\bibitem{Bao_KdV}
F.~Bao, Y.~Cao, X.~Han, and J.~Li.
\newblock Efficient particle filtering for stochastic korteweg-de vries
  equations.
\newblock {\em Stochastics and Dynamics}, 17(2):1750008, 2017.

\bibitem{BaoC2015}
F.~Bao, Y.~Cao, C.~Webster, and G.~Zhang.
\newblock A hybrid sparse-grid approach for nonlinear filtering problems based
  on adaptive-domain of the {Z}akai equation approximations.
\newblock {\em SIAM/ASA J. Uncertain. Quantif.}, 2(1):784--804, 2014.

\bibitem{Bao_BDSDE_Filter}
F.~Bao, Y.~Cao, and W.~Zhao.
\newblock A backward doubly stochastic differential equation approach for
  nonlinear filtering problems.
\newblock {\em Communications in Computational Physics}, to appear.

\bibitem{Bao_Vasilis}
F.~Bao and V.~Maroulas.
\newblock Adaptive meshfree backward sde filter.
\newblock {\em SIAM Journal on Scientific Computing}, 39(6):A2664--A2683, 2017.

\bibitem{BSDE-PIDE}
G.~Barles, R.~Buckdahn, and E.~Pardoux.
\newblock Backward stochastic differential equations and integral-partial
  differential equations.
\newblock {\em Stochastics Stochastics Rep.}, 60(1-2):57--83, 1997.

\bibitem{Bensoussan2009}
A.~Bensoussan, J.~Keppo, and S.~Sethi.
\newblock Optimal consumption and portfolio decisions with partially observed
  real prices.
\newblock {\em Math. Finance}, 19(2):215--236, 2009.

\bibitem{particle-filter-resample}
M.~Boli{\'c}, P.~Djuri{\'c}, and S~Hong.
\newblock Resampling algorithms and architectures for distributed particle
  filters.
\newblock {\em IEEE Trans. Signal Process.}, 53(7):2442--2450, 2005.

\bibitem{Chen2003}
Z.~Chen.
\newblock Bayesian filtering: From kalman filters to particle fitlers, and
  beyond.
\newblock {\em Statistics}, pages 1--69, 2003.

\bibitem{Chorin2013}
A.~Chorin, M.~Morzfeld, and X.~Tu.
\newblock A survey of implicit particle filters for data assimilation.
\newblock In {\em State-space models}, Stat. Econom. Finance, pages 63--88.
  Springer, New York, 2013.

\bibitem{MR1847785}
D.~Crisan.
\newblock Particle filters---a theoretical perspective.
\newblock In {\em Sequential {M}onte {C}arlo methods in practice}, Stat. Eng.
  Inf. Sci., pages 17--41. Springer, New York, 2001.

\bibitem{MR1895071}
D.~Crisan and A.~Doucet.
\newblock A survey of convergence results on particle filtering methods for
  practitioners.
\newblock {\em IEEE Trans. Signal Process.}, 50(3):736--746, 2002.

\bibitem{Dawson2015}
P.~Dawson, R.~Gailis, and A.~Meehan.
\newblock Detecting disease outbreaks using a combined {B}ayesian network and
  particle filter approach.
\newblock {\em J. Theoret. Biol.}, 370:171--183, 2015.

\bibitem{Djogatovi2014}
M.~Djogatovi{\'c}, M.~Stanojevi{\'c}, and N.~Mladenovi{\'c}.
\newblock A variable neighborhood search particle filter for bearings-only
  target tracking.
\newblock {\em Comput. Oper. Res.}, 52(part B):192--202, 2014.

\bibitem{Doucet2008}
A.~Doucet and A.~Johansen.
\newblock A tutorial on particle filtering and smoothing: Fifteen years later.
\newblock {\em http://www.cs.ubc.ca/~arnaud/doucet\_johansen\_tutorialPF.pdf},
  2008.

\bibitem{Elliott2013}
R.~Elliott and T.~Siu.
\newblock Option pricing and filtering with hidden {M}arkov-modulated pure-jump
  processes.
\newblock {\em Appl. Math. Finance}, 20(1):1--25, 2013.

\bibitem{RBF}
G.~F. Fasshauer.
\newblock {\em Meshfree Approximation Methods with MATLAB}.
\newblock Interdisciplinary Mathematical Sciences (Book 6). World Scientific
  Publishing Company, 2007.

\bibitem{Frey2012}
R.~Frey and T.~Schmidt.
\newblock Pricing and hedging of credit derivatives via the innovations
  approach to nonlinear filtering.
\newblock {\em Finance Stoch.}, 16(1):105--133, 2012.

\bibitem{Giessibl_2003}
F.~J. Giessibl.
\newblock Advances in atomic force microscopy.
\newblock {\em Rev. Mod. Phys.}, 75:949, 2003.

\bibitem{Gobet-Zakai}
E.~Gobet, G.~Pag{\`e}s, H.~Pham, and J.~Printems.
\newblock Discretization and simulation of the {Z}akai equation.
\newblock {\em SIAM J. Numer. Anal.}, 44(6):2505--2538 (electronic), 2006.

\bibitem{particle-filter}
N.J Gordon, D.J Salmond, and A.F.M. Smith.
\newblock Novel approach to nonlinear/non-gaussian bayesian state estimation.
\newblock {\em IEE PROCEEDING-F}, 140(2):107--113, 1993.

\bibitem{Green_2014}
M.F. Green, T.~Esat, C.~Wagner, P.~Leinen, A.~Grotsch, F.S. Tautz, and
  R.~Temirov.
\newblock Patterning a hydrogen-bonded molecular monolayer with a
  hand-controlled scanning probe microscope.
\newblock {\em Beilstein J Nanotechnol}, pages 1926--1932, 2014.

\bibitem{Hairer2011}
M.~Hairer, A.~Stuart, and J.~Voss.
\newblock Signal processing problems on function space: {B}ayesian formulation,
  stochastic {PDE}s and effective {MCMC} methods.
\newblock {\em The {O}xford handbook of nonlinear filtering}, pages 833--873,
  2011.

\bibitem{HanX2015}
X.~Han, J.~Li, and D.~Xiu.
\newblock Error analysis for numerical formulation of particle filter.
\newblock {\em Discrete Contin. Dyn. Syst. Ser. B}, 20(5):1337--1354, 2015.

\bibitem{Hastings}
W.K. Hastings.
\newblock Monte carlo sampling methods using markov chains and their
  applications.
\newblock {\em Biometrika}, 57:97?109.

\bibitem{Heller_1994}
E.J. Heller, M.F. Crommie, C.P. Lutz, and D.M. Bigler.
\newblock Scattering and absorption of surface electron waves in quantum
  corrals.
\newblock {\em Nature}, 369:464--466, 2013.

\bibitem{Kloeden_Jump}
D.~Higham and P.~Kloeden.
\newblock Numerical methods for nonlinear stochastic differential equations
  with jumps.
\newblock {\em Numerische Mathematik}, 101:101?119.

\bibitem{Hla_2014}
S.W. Hla.
\newblock Atom-by-atom assembly.
\newblock {\em Reports on Progress in Physics}, 77(5):056502, 2014.

\bibitem{HU-Zakai}
Y.~Hu, G.~Kallianpur, and J.~Xiong.
\newblock An approximation for the {Z}akai equation.
\newblock {\em Appl. Math. Optim.}, 45(1):23--44, 2002.

\bibitem{Kang_2018}
K.~Kang, V.~Maroulas, I.~Schizas, and F.~Bao.
\newblock Improved distributed particle filters for tracking in a wireless
  sensor network.
\newblock {\em Computational Statistics and Data Analysis}, 117(Supplement
  C):90 -- 108, 2018.

\bibitem{Kim2013}
D.~Kim, T.~Song, and D.~Musicki.
\newblock Highest probability data association for multi-target particle
  filtering with nonlinear measurements.
\newblock {\em IEICE TRANSACTIONS on Communications}, E96-B (1):281--290, 2013.

\bibitem{Kim_2015}
Y~Kim, K.~Motobayashi, T.~Frederiksen, H.~Ueba, and M.~Kawai.
\newblock Action spectroscopy for single-molecule reactions ? experiments and
  theory.
\newblock {\em Progress in Surface Science}, 90:85--143, 2015.

\bibitem{MR0180407}
H.~Kushner.
\newblock On the differential equations satisfied by conditional probablitity
  densities of {M}arkov processes, with applications.
\newblock {\em J. Soc. Indust. Appl. Math. Ser. A Control}, 2:106--119, 1964.

\bibitem{Lee2013}
Y.~Lee.
\newblock An asymmetric information modeling framework for ultra-high frequency
  transaction data: a nonlinear filtering approach.
\newblock In {\em State-space models}, Stat. Econom. Finance, pages 279--309.
  Springer, New York, 2013.

\bibitem{Little-Jones2013}
M.~Little and N.~Jones.
\newblock Signal processing for molecular and cellular biological physics: an
  emerging field.
\newblock {\em Philos. Trans. R. Soc. Lond. Ser. A Math. Phys. Eng. Sci.},
  371(1984):20110546, 18, 2013.

\bibitem{Peter_book}
P.~Maksymovych, S.~V. Kalinin, and A.~Gruverman.
\newblock {\em Scanning probe microscopy of functional materials nanoscale
  imaging and spectroscopy}.
\newblock Springer, 2011.

\bibitem{Monte-Carlo}
N.~Metropolis and S.~Ulam.
\newblock The monte carlo method.
\newblock {\em Journal of the American Statistical Association}, 44:335--341.

\bibitem{Jump-finance}
T.~Meyer-Brandis and F.~Proske.
\newblock Explicit solution of a non-linear filtering problem for {L}\'evy
  processes with application to finance.
\newblock {\em Appl. Math. Optim.}, 50(2):119--134, 2004.

\bibitem{Morgenstern_2013}
N.~Morgenstern, N.~Lorente, and K.-H Rieder.
\newblock Controlled manipulation of single atoms and small molecules using the
  scanning tunnelling microscope.
\newblock {\em Phys. Status Solidi B}, 250:1671?1751, 2013.

\bibitem{Pardoux1990}
{\'E}.~Pardoux and S.~G. Peng.
\newblock Adapted solution of a backward stochastic differential equation.
\newblock {\em Systems Control Lett.}, 14(1):55--61, 1990.

\bibitem{Two_sided}
{\'E}.~Pardoux and P.~Protter.
\newblock A two-sided stochastic integral and its calculus.
\newblock {\em Probab. Theory Related Fields}, 76(1):15--49, 1987.

\bibitem{APF}
M.. Pitt and N.~Shephard.
\newblock Filtering via simulation: auxiliary particle filters.
\newblock {\em J. Amer. Statist. Assoc.}, 94(446):590--599, 1999.

\bibitem{Platen_Jump}
E.~Platen and N.~Bruti-Liberati.
\newblock {\em Numerical Solution of Stochastic Differential Equations with
  Jumps in Finance}.
\newblock Springer, 2010.

\bibitem{Popa-jump}
S.~Popa and S.~S. Sritharan.
\newblock Nonlinear filtering of {I}t\^o-{L}\'evy stochastic differential
  equations with continuous observations.
\newblock {\em Commun. Stoch. Anal.}, 3(3):313--330, 2009.

\bibitem{Duan}
H.~Qiao and J.~Duan.
\newblock Nonlinear filtering of stochastic dynamical systems with lévy noises.
\newblock {\em Advances in Applied Probability}, 47(3):902?918, 2015.

\bibitem{Singh2013}
S.~Singh, S.~Digumarthy, A.~Back, J.~Shepard, and Kalra M.
\newblock Radiation dose reduction for chest ct with non-linear adaptive
  filters.
\newblock {\em Acta Radiologica}, 55(2):169--174, 2013.

\bibitem{Stepanov2015}
O.~A. Stepanov, D.~A. Koshaev, and A.~V. Motorin.
\newblock Identification of gravity anomaly model parameters in airborne
  gravimetry problems using nonlinear filtering methods.
\newblock {\em Gyroscopy and Navigation}, 6(4):318--323, 2015.

\bibitem{stratonovich1960}
R.~Stratonovich.
\newblock Conditional markov processes.
\newblock {\em Theoretical Probability and Its Applications}, 5:156--178, 1960.

\bibitem{Duan_F-P}
X.~Sun and J.~Duan.
\newblock Fokker-{P}lanck equations for nonlinear dynamical systems driven by
  non-{G}aussian {L}\'evy processes.
\newblock {\em J. Math. Phys.}, 53(7):072701, 10, 2012.

\bibitem{Ternes_2008}
M.~Ternes, C.P. Lutz, F.J. Hirjibehedin, F.J. Giessibl, and A.J. Heinrich.
\newblock The force needed to move an atom on a surface.
\newblock {\em Science}, 319:1066--1069, 2008.

\bibitem{Walter-Error}
G.~Walter.
\newblock Norm estimate for inverses of vandermonde matrices.
\newblock {\em Numerische Mathematik}, 23(4):337?347, 1975.

\bibitem{Wea-PF}
J.~Weare.
\newblock Particle filtering with path sampling and an application to a bimodal
  ocean current model.
\newblock {\em J. Comput. Phys.}, 228(12):4312--4331, 2009.

\bibitem{zakai}
M.~Zakai.
\newblock On the optimal filtering of diffusion processes.
\newblock {\em Ztschr. Wahrscheinlichkeitstheor und verw. Geb.}, 11:230--243,
  1969.

\bibitem{zhang}
H.~Zhang and D.~Laneuville.
\newblock Grid based solution of zakai equation with adaptive local refinement
  for bearing-only tracking.
\newblock {\em IEEE Aerospace Conference}, 2008.

\bibitem{Karniadakis_Temp}
M.~Zheng and G.~Karniadakis.
\newblock Numerical methods for spdes with tempered stable processes.
\newblock {\em SIAM Journal of Scientific Computing}, 37(3):A1197?A1217, 2015.

\end{thebibliography}

\def\cprime{$'$} \def\cprime{$'$}

\end{document}